\input amstex
\documentstyle{amsppt}
\NoBlackBoxes
\mag1100
\let\aro=@ 
\newcount\numcount
\def\numerote{\global\advance\numcount by 1 \the\numcount}
\def\lastnum[#1]{{\advance \numcount by #1 \the\numcount}}
\loadbold
%\input psbox 
%\def \scaledpicture #1 by #2 (#3 scaled #4) 
%{\dimen0=#1 \dimen1=#2 \divide\dimen0 by 1000\multiply\dimen0 by #4 
%\divide\dimen1 by 1000\multiply\dimen1 by #4 
%$$\psboxto(\dimen0;\dimen1){#3.ps}$$} 
\topmatter
\title 
Vari\'et\'es CR polaris\'ees et $\boldkey G$-polaris\'ees, partie I
\endtitle 
\rightheadtext{Vari\'et\'es CR polaris\'ees I}
\author
Laurent Meersseman
\endauthor
\date 8 f\'evrier 2013\enddate 
\address 
\noindent{Institut de Math\'ematiques de Bourgogne, B\^atiment Mirande, B.P. 47870, 21078 DIJON Cedex, FRAN\-CE}
\endaddress
\curraddr
Centre de Recerca Matem\`atica, Edifici C, Campus de Bellaterra, 08193 BELLATERRA, ESPA\~ NA
\endcurraddr
\email laurent.meersseman\@u-bourgogne.fr {\it and} Lmeersseman\@crm.cat\endemail
\dedicatory
\`A la m\'emoire de Marco Brunella.
\enddedicatory
\abstract
Polarized and $G$-polarized CR manifolds are smooth manifolds endowed with a double structure: a real foliation $\Cal F$ (given by the action of a Lie group $G$ in the $G$-polarized case) and a transverse CR distribution $(E,J)$. Polarized means that $(E,J)$ is roughly speaking invariant by $\Cal F$. Both structures are therefore linked up. The interplay between them gives to polarized CR-manifolds a very rich geometry.

In this paper, we study the properties of polarized and $G$-polarized manifolds, putting special emphasis on the existence of a local moduli space.
\endabstract
\thanks
Ce travail a une longue histoire. Il a d\'ebut\'e en 2008-2009 au cours d'une d\'el\'egation au CNRS d'un an au Pacific Institute for Mathematical Sciences, sur le campus de la University of British Columbia, Vancouver, BC. Il s'est poursuivi \`a Dijon ensuite. Mais il n'a vraiment pris forme que deux ans plus tard, au Centre de Recerca Matem\`atica de Bellaterra, comme partie principale du projet Marie Curie 27707 DEFFOL.
\smallskip
Je voudrais remercier d'une part le PIMS et le d\'epartement de math\'ematiques de UBC, d'autre part le CRM pour leur hospitalit\'e. Merci \'egalement au CNRS qui m'a donn\'e l'occasion de cette d\'el\'egation. Merci enfin \`a Marcel Nicolau et \`a Graham Smith pour de nombreuses et fructueuses discussions au cours de l'\'elaboration de ce travail.
\smallskip
Ce travail a \'et\'e financ\'e principalement par la bourse Marie Curie 27707 DEFFOL de la communaut\'e europ\'eenne et, avant cela, partiellement par le projet COMPLEXE (ANR-08-JCJC-0130-01) de l'Agence Nationale de la Recherche. 
\endthanks
\endtopmatter
\document
\vskip1cm
\head 
{\bf Introduction.}
\endhead

Les vari\'et\'es CR polaris\'ees et les vari\'et\'es CR $G$-polaris\'ees sont des vari\'et\'es lisses munies d'une double structure : un feuilletage r\'eel $\Cal F$ et une distribution CR transverse $(E,J)$ (par CR nous entendons CR int\'egrable). On suppose la distribution polaris\'ee par le feuilletage, autrement dit poss\'edant une propri\'et\'e d'invariance par le feuilletage. Ceci lie intimement la composante r\'eelle et la composante complexe de ces objets mixtes $C^\infty$/holomorphe. L'interaction r\'eel/complexe leur donne une g\'eom\'etrie tr\`es riche et des propri\'et\'es de d\'eformations remarquables.
\medskip
Dans le cas $G$-polaris\'e, $\Cal F$ est donn\'e par l'action d'un groupe de Lie r\'eel $G$ et la propri\'et\'e de polarisation signifie que $G$ pr\'eserve $(E,J)$. On obtient ainsi une notion a priori \`a mi-chemin entre vari\'et\'e complexe (lorsque la dimension de $G$ est nulle) et vari\'et\'e r\'eelle homog\`ene (lorsqu'elle est maximale). Toutefois, une analyse plus approfondie de la situation montre que ces objets sont beaucoup plus proches de la g\'eom\'etrie complexe que de la g\'eom\'etrie r\'eelle. D'une part, si $G$ est compact de dimension paire, alors une vari\'et\'e CR $G$-polaris\'ee $X$ est complexe avec $\Cal F$ holomorphe (ce qui forme un cas peu int\'eressant). D'autre part, si $G$ est ab\'elien, le produit $X\times G$ poss\`ede un atlas de vari\'et\'e complexe pour lequel les translations sur les fibres de $X\times G\to X$ sont holomorphes et ces fibres sont totalement r\'eelles. Cette propri\'et\'e fondamentale, qui signifie que le feuilletage r\'eel $\Cal F$ peut \^etre complexifi\'e au sens de \cite{H-S}, rapproche les vari\'et\'es $\Bbb R^m$-polaris\'ees des vari\'et\'es sasakiennes dont, rappelons-le, le produit avec $\Bbb R^{>0}$ admet une structure complexe k\"ahl\'erienne invariante par dilatations ; mais sans aspect m\'etrique ni forme de contact. Mais elle les relie \'egalement \`a la g\'eom\'etrie complexe non k\"ahl\'erienne. En particulier, on obtient des exemples triviaux de structures $\Bbb R^m$-polaris\'ees en consid\'erant une somme de Whitney de fibr\'es unitaires associ\'ee \`a une somme de Whitney de $m$ fibr\'es en droites sur une vari\'et\'e compacte complexe. Et, ce qui est beaucoup plus int\'eressant, on obtient des exemples non-triviaux en consid\'erant des quotients de certaines vari\'et\'es LVM (cf. \cite{LdM-V}, \cite{Me1} et \cite{M-V} pour une pr\'esentation de cette classe de vari\'et\'es non k\"ahl\'eriennes).
\medskip
Dans le cas plus g\'en\'eral des vari\'et\'es polaris\'ees, il n'y a pas d'action de groupe. La distribution CR transverse est invariante par holonomie, au sens o\`u elle induit une structure transverse holomorphe pour le feuilletage. Cette propri\'et\'e de polarisation entra\^{\i}ne une propri\'et\'e de rigidit\'e : ces structures poss\`ede un espace de modules local\footnote{On parle ici d'espace de modules \`a type constant, pour reprendre la terminologie introduite en section 2, c'est-\`a-dire qu'on fixe la distribution et qu'on d\'eforme uniquement l'op\'erateur CR sur cette distribution.} (au sens de Kuranishi) de dimension finie. La finitude est un quasi-miracle pour un espace de structures CR (on consultera \`a ce propos l'introduction de \cite{Me2}). Elle refl\`ete bien s\^ur l'existence d'un feuilletage transversalement holomorphe associ\'e \`a la structure polaris\'ee ; cependant, ce r\'esultat ne se ram\`ene pas aux \'enonc\'es classiques de \cite{EK-N} et \cite{Gi} (qui d'ailleurs n\'ecessitent des hypoth\`eses suppl\'ementaires) et sa preuve ne les utilise pas. De surcro\^{\i}t, cet espace est {\it g\'eom\'etrique}. Les submersions \`a base complexe fournissent le meilleur exemple de ce ph\'enom\`ene. Ici la distribution CR donne naissance \`a un feuilletage Levi-plat. En tant que structure polaris\'ee, leur espace de modules local s'identifie \`a l'espace de Kuranishi de la base et ne d\'epend donc pas de la repr\'esentation d\'efinissant la suspension. En tant que structure CR, leur espace de modules\footnote{Si tant est qu'on arrive \`a donner un sens pr\'ecis \`a ce terme.} est g\'en\'eralement de dimension infinie, mais des variations infimes de la repr\'esentation peuvent induire un espace de modules de dimension finie, et ce {\it sans que la g\'eom\'etrie du feuilletage Levi-plat ne soit affect\'ee}. Ainsi, il n'y a pas moyen de lire la finitude de l'espace des modules sur la g\'eom\'etrie du feuilletage Levi-plat.
\medskip
Il y a cependant un prix \`a payer pour obtenir cette finitude. Il faut identifier deux structures non pas modulo CR isomorphisme, mais modulo une relation d'\'equivalence plus grossi\`ere, qui, grosso modo, signifie que deux structures \'equivalentes sont CR isomorphes {\it \`a l'ordre un}. Agr\'eablement, sous une condition m\'etrique classique, on identifie les structures si elles induisent la m\^eme structure transverse holomorphe sur $\Cal F$. Enfin pour les submersions \`a base complexe, l'identification se fait modulo CR isomorphisme.
\medskip
L'article comprend trois parties. La premi\`ere contient une description infinit\'e\-si\-male des structures CR polaris\'ees, les premiers exemples et tous les pr\'eliminaires n\'ecessaires \`a la construction d'un espace de modules \`a la Kuranishi. La deuxi\`eme construit cet espace de modules, culminant avec les th\'eor\`emes 10.1, 13.1 et leurs corollaires ; et finit par la description du cas Levi-plat, et de celui des submersions \`a base complexe (th\'eor\`eme 15.2). La troisi\`eme d\'efinit et \'etudie les vari\'et\'es CR $G$-polaris\'ees. On montre en section 17 que, dans le cas ab\'elien, le produit par $G$ d'une telle vari\'et\'e, et plus g\'en\'eralement tout $G$-fibr\'e plat au-dessus d'elle, admet une structure complexe invariante par les translations du groupe et pour laquelle les fibres sont totalement r\' eelles. Et on prouve en section 18 que, dans le cas $G$ compact de dimension paire, une vari\'et\'e CR $G$-polaris\'ee est une vari\'et\'e complexe munie d'un feuilletage holomorphe. On finit l'article par l'exemple non-trivial des structures $\Bbb R^m$-polaris\'ees issues des vari\'et\'es LVM.
\medskip
Sur le plan technique, les preuves des th\'eor\`emes 10.1 et 13.1 suivent le sch\'ema classique de \cite{Ku3} (espace de modules local de structures complexes) ou de \cite{Do-K} (espace de modules local de connexions anti autoduales). Il faut toutefois noter que l'op\'erateur que nous utilisons (et dont le noyau est tangent \`a l'espace de modules construit) {\it n'est pas un op\'erateur diff\'erentiel}, mais simplement un op\'erateur lin\'eaire entre Banach. En cons\'equence, il n'y a pas de r\'esolution elliptique ni de faisceau de champs de vecteurs naturellement associ\'e \`a notre situation. Cette diff\'erence d'outils ne modifie pas la structure de l'argument (au lieu de montrer qu'un op\'erateur diff\'erentiel est elliptique, on montre qu'un op\'erateur lin\'eaire est Fredholm), mais s'av\`ere cruciale pour gagner en flexibilit\'e et obtenir les th\'eor\`emes dans toute leur g\'en\'eralit\'e (cf. le d\'ebut du paragraphe 11).
\medskip
Un deuxi\`eme article, incluant une \'etude plus pouss\'ee des vari\'et\'es $G$-polaris\'ees, ainsi que l'analyse des d\'eformations \`a polarisation constante (c'est-\`a-dire sans fixer le fibr\'e tangent de la structure CR), est pr\'evu.

\vfill\eject

\head
I. G\'en\'eralit\'es.
\endhead
\vskip1cm
\subhead
{\bf 1. Structures presque CR polaris\'ees}
\endsubhead

Soit $X$ une vari\'et\'e r\'eelle lisse (i.e. $C^\infty$), compacte, connexe de dimension $n$. Soit $TX$ son fibr\'e tangent et soit 
$$
T_{{\Bbb C}}X=TX\otimes_{\Bbb R}\Bbb C
$$
son fibr\'e tangent complexifi\'e.

\remark{Remarque}
Nous insistons sur le fait que la dimension de $X$ est {\it arbitraire}. En particulier, $X$ peut \^etre {\it de dimension impaire}.
\endremark
\medskip

Soit $(E,J)$ une structure presque CR de classe $C^{\infty}$ de $X$.

\definition{D\'efinition}
Une {\it polarisation} of $(E,J)$ est un suppl\'ementaire $N$ de $E$ dans $TX$. 
Une structure presque CR munie d'une polarisation est dite {\it polaris\'ee par} $N$ ou plus simplement {\it polaris\'ee} lorsqu'il n'y a pas d'ambigu\"{\i}t\'e.
\enddefinition

Soit $E_{\Bbb C}$ le complexifi\'e du sous-fibr\'e $E$. Soit $J_{\Bbb C}$ l'extension de l'op\'erateur $J$ \`a $E_{\Bbb C}$. On a une d\'ecomposition
$$
E_{\Bbb C}=E^{0,1}\oplus E^{1,0}
\leqno (1.\numerote)$$
o\`u $E^{0,1}$, respectivement $E^{1,0}$, est le sous-fibr\'e propre associ\'e \`a la valeur propre $-i$, respectivement $i$. En tenant compte de la polarisation, ceci donne
$$
\eqalign{
T_{\Bbb C}X=&E^{0,1}\oplus E^{1,0}\oplus N\oplus iN\cr
=&(E^{0,1}\oplus iN)\oplus i(\overline{E^{0,1}\oplus iN})\cr
:=&T\oplus i\overline{T}
}
\leqno (1.\numerote)
$$

Soit $G\to X$ le fibr\'e r\'eel lisse dont la fibre en $x\in X$ est la grassmannienne des $n$-plans {\it r\'eels} de l'espace vectoriel complexe $(T_{\Bbb C}X)_{x}$. Dans (1.\lastnum[0]), le fibr\'e $T$ est une section lisse de $G$.
\medskip
Remarquons que $T$ d\'efinit compl\`etement la structure presque CR polaris\'ee associ\'ee $(E,J,N)$. En effet, on a
$$
E^{0,1}=T\cap iT\qquad\qquad E^{1,0}=\overline{E^{0,1}}
\leqno (1.\numerote)
$$
ce qui d\'efinit $E_\Bbb C$ via (1.\lastnum[-2]) et donc $E$. Mais on d\'efinit alors $J_\Bbb C$ comme la multiplication par $i$ sur $E^{1,0}$ et par $-i$ sur $E^{0,1}$. En restriction \`a $E$, l'op\'erateur $J_\Bbb C$ devient $J$. Enfin, un calcul imm\'ediat montre que
$$
N=i\overline{T}\cap TX.
\leqno (1.\numerote)
$$
On vient de montrer

\proclaim{Proposition 1.1}
Les structures presque CR polaris\'ees (lisses) sont en bijection avec les sections (lisses) $T$ de $G$ telles que 
$$
\eqalign{
(i) \qquad &T\oplus i\overline{T}=T_{\Bbb C}X\cr
(ii) \qquad &T\cap i{T} \text{ est un sous-fibr\'e complexe de }T_{\Bbb C}X\cr  
(iii) \qquad &TX\cap i\overline{T} \text{ est un sous-fibr\'e r\'eel de }TX\cr
(iv) \qquad &T=E^{0,1}\oplus iN
}
$$
\endproclaim

\vskip.5cm
\subhead
{\bf 2. Structures proches}
\endsubhead
\numcount=0

On munit $\Sigma (G)$, l'espace des {\it distributions} de $G$, c'est-\`a-dire des sections lisses de $G$, d'une norme Sobolev $L^2_l$ (pour $l>0$ un entier donn\'e) de la mani\`ere classique (cf. \cite{Ku2, \S IX}). 
\medskip
Soit $T$ une structure presque CR polaris\'ee. Soit $T'$ une distribution proche de $T$ en norme $l$. On peut voir $T'$ comme un graphe de $T$ dans $i\overline{T}$, i.e. il existe 
$$
\omega\in \text{Hom}_{\Bbb R}(T,i\overline{T})
\leqno (2.\numerote)
$$
telle que
$$
T'=\{v-\omega(v)\qquad\vert\qquad v\in T\}
\leqno (2.\numerote)
$$
Cependant, une distribution (2.\lastnum[0]) ou, de fa\c con \'equivalente, un morphisme (2.\lastnum[-1]), ne code pas forc\'ement une structure presque CR polaris\'ee. Le sous-fibr\'e (2.\lastnum[0]) doit v\'erifier les conditions (i), (ii), (iii) et (iv) de la proposition 1.1.
\medskip
En fait, nous allons consid\'erer deux cas particuliers distincts de structures pro\-ches, pour lesquelles on peut caract\'eriser agr\'eablement le morphisme (2.\lastnum[-1]) associ\'e.

\medskip
Soit $T$ une structure presque CR polaris\'ee. Soit $T'$ une structure presque CR polaris\'ee proche de $T$ en norme $l$. Soit $\omega$ le morphisme associ\'e \`a $T'$ via (2.\lastnum[0]).
Appelons $E'$, $E'_\Bbb C$, $(E^{0,1})'$ et $N'$ les sous-fibr\'es associ\'es \`a $T'$.

\definition{D\'efinitions}
On dira que $T'$ est une {\it d\'eformation de $T$ \`a type constant} si l'on a $E'\equiv E$ et $N'\equiv N$.
\medskip
On dira que $T'$ est une {\it d\'eformation de $T$ \`a polarisation constante} si l'on a simplement $N'\equiv N$.
\enddefinition

Ainsi dans les deux cas, on laisse la polarisation fixe. Mais dans le premier, on se contente de modifier l'op\'erateur presque complexe $J$ le long de $E$, qui, lui, reste fixe ; tandis que dans le second, on modifie $E$ et $J$.  
\medskip
Pla\c cons-nous \`a type constant. Alors le morphisme $\omega$ donn\'e par (2.\lastnum[-1]) est nul sur $N$, qui ne varie pas. Qui plus est, $(E^{0,1})'$ est un sous-fibr\'e complexe de $E'_\Bbb C$, si bien qu'il existe
$$
\omega\in \text{Hom}_{\Bbb C}(E^{0,1},E^{1,0})
\leqno (2.\numerote)
$$
tel que l'\'egalit\'e (2.\lastnum[-1]) se r\'e\'ecrive
$$
T'=(E^{0,1})'\oplus iN=\{v-\omega(v)\quad\vert\quad v\in E^{0,1}\}\oplus iN.
\leqno (2.\numerote)
$$
Et, ce qui est beaucoup plus int\'eressant, on obtient cette fois une caract\'erisation compl\`ete des d\'eformations \`a type constant.

\proclaim{Proposition 2.1}
L'espace des d\'eformations de $T$ \`a type constant est en bijection avec l'espace des morphismes complexes de $E^{0,1}$ dans $E^{1,0}$ via l'identit\'e {\rm (2.\lastnum[0])}.
\endproclaim

\demo{Preuve}
Tout d\'ecoule de la formule (2.\lastnum[0]) coupl\'ee \`a la proposition 1.1.
$\square$
\enddemo

Fixons $E$ et $N$. Soit
$\Cal E_{tc}$ l'ensemble des structures presque CR lisses polaris\'ees $(E,J,N)$. Il s'agit d'un sous-ensemble de $\Sigma (G)$ que nous munissons de la topologie induite.
\medskip
Compl\'etons $\Sigma (G)$ en $\Sigma_l (G)$ pour la norme $l$ ; l'espace $\Cal E_{tc}$ se compl\`ete en $\Cal E_{tc}^l$. La proposition 2.1 a l'int\'eressant corollaire suivant.

\proclaim{Corollaire 2.2}
L'espace $\Cal E_{tc}^l$ est une vari\'et\'e de Banach $\Bbb C$-analytique.
\endproclaim

\demo{Preuve}
On observe que la proposition 2.1 donne des cartes de vari\'et\'e de Banach $\Bbb C$-analytique. En effet, la bijection annonc\'ee est en fait un hom\'eomorphisme pour la topologie $L^2_l$. Autrement dit, $\Cal E^l_{tc}$ a pour carte locale en 
$$
T=E^{0,1}\oplus iN
$$
un voisinage de $0$ dans l'espace de Banach Hom$^l_\Bbb C(E^{0,1},E^{1,0})$ pour la norme $l$.
\medskip
Par ailleurs, les changements de cartes de cet atlas sont simplement des changements de cartes de $G$. Qui plus est, ces changements de cartes ne touchent pas \`a la polarisation $N$ et peuvent donc \^etre choisis $\Bbb C$-analytiques. 
$\square$
\enddemo
Nous consid\'erons maintenant les d\'eformations \`a polarisation constante. Comme pr\'ec\'edemment, le morphisme $\omega$ donn\'e par (2.\lastnum[-3]) est nul sur $N$, qui ne varie pas. Ceci permet de prendre
$$
\omega\in \text{Hom}_\Bbb R (E^{0,1}, i\overline{T})
$$
mais on ne peut plus cette fois supposer $\omega$ complexe \`a valeurs dans $E^{1,0}$. D\'e\-com\-po\-sons
$$
\left\{\eqalign{
&\omega=\omega_{1,0}+\omega_N\cr
&\text{avec }\omega_{1,0}\in \text{Hom}_{\Bbb R}(E^{0,1},E^{1,0})\text{ et }\omega_N\in\text{Hom}_{\Bbb R}(E^{0,1},N) 
}
\right .
\leqno (2.\numerote)
$$

Soit $w\in (E^{0,1})'$. Par (2.\lastnum[-2]), il existe $v\in E^{0,1}$ et $n\in N$ tels que
$$
w=v\oplus\omega_{1,0}(-v)\oplus\omega_N(-v)\oplus in \in E^{0,1}\oplus E^{1,0}\oplus N\oplus iN
$$
D\`es lors,
$$
iw=iv\oplus i\omega_{1,0}(-v)\oplus (-n)\oplus i\omega_N(-v)\in E^{0,1}\oplus E^{1,0}\oplus N\oplus iN=T_\Bbb C X
\leqno (2.\numerote)
$$
Mais $iw\in (E^{0,1})'$ et il existe donc $u\in E^{0,1}$ et $p\in N$ tels que
$$
iw=u\oplus \omega_{1,0}(-u)\oplus \omega_N(-u)\oplus ip\in E^{0,1}\oplus E^{1,0}\oplus N\oplus iN
\leqno (2.\numerote)
 $$
En comparant (2.\lastnum[0]) et (2.\lastnum[-1]), et en notant que la projection
$$
w\in (E')^{0,1}\longmapsto v\in E^{0,1}
$$
est un isomorphisme, on conclut que $\omega_{1,0}$ est en fait une forme complexe ; et que
$$
\omega_{N_\Bbb C}(v):=\omega_N(v)+i\omega_N(iv)\quad\text{ pour }v\in E^{0,1}
\leqno (2.\numerote)
$$
est elle aussi une forme complexe.
\medskip
En somme, nous voyons qu'on peut r\'e\'ecrire (2.\lastnum[-6]) comme
$$
T'=(E^{0,1})'\oplus iN=\{v-\omega_\Bbb C(v)\quad\vert\quad v\in E^{0,1}\}\oplus iN.
\leqno (2.\numerote)
$$
o\`u
$$
\left\{\eqalign{
&\omega_{\Bbb C}:=\omega_{1,0}+\omega_{N_\Bbb C}\cr
&\text{avec }\omega_{1,0}\in \text{Hom}_{\Bbb C}(E^{0,1},E^{1,0})\text{ et }\omega_{N_\Bbb C}\in\text{Hom}_{\Bbb C}(E^{0,1},N_{\Bbb C}) 
}
\right .
\leqno (2.\numerote)
$$
Ce sont les formules \'equivalentes, pour une d\'eformation \`a polarisation constante, aux formules (2.\lastnum[-7]) et (2.\lastnum[-6]).
\medskip
Nous pouvons maintenant \'enoncer

\proclaim{Proposition 2.3}
L'espace des d\'eformations de $T$ \`a polarisation constante est en bijection avec l'espace des morphismes complexes de $E^{0,1}$ dans $E^{1,0}\oplus N_\Bbb C$ via l'identit\'e {\rm (2.\lastnum[-1])}.
\endproclaim

\demo{Preuve}
Tout d\'ecoule de la formule (2.\lastnum[-1]) coupl\'ee \`a la proposition 1.1.
$\square$
\enddemo

Fixons $N$ et d\'efinissons
$\Cal E_{pc}$ comme l'ensemble des structures presque CR lisses polaris\'ees $(E,J,N)$. Compl\'etons-le en  $\Cal E_{pc}^l$ pour la norme $l$. On montre de la mani\`ere que le corollaire 2.2 le

\proclaim{Corollaire 2.4}
L'espace $\Cal E_{pc}^l$ est une vari\'et\'e de Banach $\Bbb C$-analytique.
\endproclaim

\vskip.5cm
\subhead
{\bf 3. Structures et vari\'et\'es CR polaris\'ees}
\endsubhead
\numcount=0

Soit $T$ une structure presque CR polaris\'ee. On note $T^\Bbb C$ la distribution de $T_\Bbb C X$ \'egale en $x\in X$ au plus petit sous-espace vectoriel complexe de $(T_\Bbb C X)_x$ contenant le sous-espace r\'eel $T_x$. Et on note $T_\Bbb C$ la distribution de $T_\Bbb C X$ \'egale en $x\in X$ au plus grand sous-espace vectoriel complexe de $(T_\Bbb C X)_x$ contenu dans le sous-espace r\'eel $T_x$.
\medskip
Il est imm\'ediat que ces deux distributions sont en fait des sous-fibr\'es complexes de $T_\Bbb C X$ et qu'on a 
$$
T^\Bbb C=E^{0,1}\oplus N_\Bbb C \qquad\text{ et }\qquad T_\Bbb C=E^{0,1}
\leqno (3.\numerote)
$$

\definition{D\'efinition}
On dit que $T$ est une structure presque CR polaris\'ee {\it compl\`etement int\'egrable}, ou, plus bri\`evement, une {\it structure CR polaris\'ee} si $T^\Bbb C$ et $T_\Bbb C$ sont involutifs, i.e. si
$$
[T^\Bbb C,T^\Bbb C]\subset T^\Bbb C \qquad\text{ et }\qquad [T_\Bbb C,T_\Bbb C]\subset T_\Bbb C
$$
\enddefinition
\definition{D\'efinition}
Une vari\'et\'e r\'eelle $X$ munie d'une structure CR polaris\'ee sera appel\'ee {\it vari\'et\'e CR polaris\'ee}.
\enddefinition

Compte-tenu de (3.\lastnum[0]), l'involutivit\'e de $T_\Bbb C$ signifie exactement l'int\'egrabilit\'e de la structure presque CR $(E,J)$. La compl\`ete int\'egrabilit\'e est donc, comme son nom le sugg\`ere, plus forte que l'int\'egrabilit\'e.
\medskip
Observons que, $T$ \'etant une distribution {\it r\'eelle} de $T_\Bbb C X$, elle ne peut \^etre involutive. Le mieux que l'on puisse demander est
$$
[T,T]\subset T^\Bbb C
$$
qui entra\^{\i}ne l'involutivit\'e de $T^\Bbb C$. C'est la condition suppl\'ementaire contenue dans la compl\`ete int\'egrabilit\'e par rapport \`a l'int\'egrabilit\'e classique.
\medskip
Convenons de dire que $(E,J,N)$ est {\it transversalement int\'egrable} si elle v\'erifie uniquement cette deuxi\`eme condition, \`a savoir l'involutivit\'e de $T^\Bbb C$. Voici le sens de l'int\'egrabilit\'e transverse.

\proclaim{Proposition 3.1}
Soit $T$ une structure presque CR polaris\'ee. Supposons $T$ trans\-versalement int\'egrable. Alors la polarisation $N$ est tangente \`a un feuilletage trans\-versalement holomorphe $\Cal F$ dont le fibr\'e normal $N\Cal F$ est $\Bbb C$-isomorphe au fibr\'e complexe $(E,J)$.
\endproclaim

\demo{Preuve}
Posons $F=T^\Bbb C$. On remarque que $F$ est involutif et v\'erifie, d'apr\`es (3.\lastnum[0])
$$
F+\overline{F}=T_\Bbb CX\qquad\text{ et }\qquad F\cap \overline{F}=N
$$
Par \cite{Nir}, on a imm\'ediatement que $N$ est le fibr\'e tangent \`a un feuilletage transversalement holomorphe $\Cal F$ dont la structure transverse est donn\'ee ainsi :
la projection naturelle
$$
\pi \ :\ TX\longrightarrow N\Cal F:=TX/T\Cal F
\leqno (3.\numerote)
$$
est un isomorphisme entre $E$ et le fibr\'e normal $N\Cal F$, qui permet de pousser la structure presque CR de $E$ en une structure CR de $N\Cal F$. Plus pr\'ecis\'ement, on \'etend $\pi$ en
une projection de $T_\Bbb CX$ \`a valeurs dans $N_\Bbb C\Cal F:=N\Cal F\otimes \Bbb C$
et on obtient une d\'ecomposition 
$$
N_\Bbb C\Cal F=N\Cal F^{0,1}\oplus N\Cal F^{1,0}
\leqno (3.\numerote)
$$
en posant
$$
N\Cal F^{0,1}=\pi_*E^{0,1}\qquad\text{ et }\qquad N\Cal F^{1,0}=\pi_*E^{1,0}.
\leqno (3.\numerote)
$$
Observons que l'involutivit\'e de $T^\Bbb C$ entra\^{\i}ne, par passage au quotient, l'involutivit\'e de $N\Cal F^{0,1}$. Ce qui ach\`eve la preuve.
$\square$
\enddemo

Ainsi, \`a toute structure presque CR polaris\'ee transversalement int\'egrable (et donc a fortiori \`a toute structure CR polaris\'ee), on associe un feuilletage transversalement holomorphe. Le lecteur non averti pourra se reporter \`a \cite{Ni} pour plus de d\'etails sur ces feuilletages.
\medskip
R\'eciproquement, \'etant donn\'e un feuilletage transversalement holomorphe $\Cal F$ de fibr\'e tangent $N$ et de fibr\'e normal $N\Cal F$, on peut associer une structure presque CR polaris\'ee par $N$ transversalement int\'egrable. En effet, comme $\Cal F$ est transversalement holomorphe, il existe une structure complexe sur $N\Cal F$. Autrement dit, on a une d\'ecomposition (3.\lastnum[-1]) du fibr\'e normal complexifi\'e
qui v\'erifie de plus
$$
[N\Cal F^{0,1},N\Cal F^{0,1}]\subset N\Cal F^{0,1}.
\leqno (3.\numerote)
$$
Dans (3.\lastnum[0]), le crochet utilis\'e est le crochet "quotient" sur $N\Cal F$. Autrement dit, (3.\lastnum[0]) signifie que, pour tous champs locaux $\xi$ et $\eta$ de $T_\Bbb CX$ dont la projection (3.\lastnum[-3]) appartient \`a $N\Cal F^{0,1}$, on a
$$
[\pi_*\xi,\pi_*\eta]:=\pi_*[\xi,\eta]\in N\Cal F^{0,1}.
\leqno (3.\numerote)
$$
Il suffit de choisir un sous-fibr\'e $E$ de $TX$ qui soit une r\'ealisation de $N\Cal F$.  La projection (3.\lastnum[-4]) induit un isomorphisme entre $E$ et $N\Cal F$ si bien que la d\'ecomposition (3.\lastnum[-3]) induit une d\'ecomposition 
$$
E_\Bbb C=E^{0,1}\oplus E^{1,0}
\leqno (3.\numerote)
$$
et l'involutivit\'e (3.\lastnum[-1]) implique
$$
[E^{0,1},E^{0,1}]\subset E^{0,1} \mod N_\Bbb C.
\leqno (3.\numerote)
$$
Maintenant, (3.\lastnum[0]) entra\^{\i}ne l'int\'egrabilit\'e transverse de $E$, en prenant en compte (3.\lastnum[-7]) et le fait que $N_\Bbb C$ est involutif, puisque tangent \`a un feuilletage.
\medskip
Le choix d'une r\'ealisation $E$ n'est cependant pas unique. On peut bien entendu rigidifier les choses en munissant $X$ d'une m\'etrique riemannienne et en prenant pour $E$ le suppl\'ementaire orthogonal de $N$ dans $TX$, mais {\it il n'existe pas de choix canonique de} $E$.
\medskip
Pour r\'esumer ce que nous venons de dire, soit $g$ une m\'etrique riemannienne sur $X$. D\'efinissons
$$
\Cal E_g=\{(E,J,N)\text{ structure presque CR polaris\'ee lisse }\quad\vert\quad E\perp N\}
\leqno (3.\numerote)
$$

C'est un sous-ensemble de $\Cal E_{pc}$. On vient de prouver

\proclaim{Corollaire 3.2}
Fixons $E$ et $N$. Alors il existe une bijection entre le sous-espace des \'el\'ements transversalement int\'egrables de $\Cal E_{tc}$  et l'ensemble des feuilletages trans\-ver\-sa\-le\-ment holomorphes tangents \`a $N$.

De plus, il existe une bijection entre le sous-ensemble des \'el\'ements transversalement int\'egrables de $\Cal E_g$ et l'ensemble des feuilletages transversalement holomorphes sur $X$.

\endproclaim

Ainsi, une fois choisie une m\'etrique riemannienne sur $X$, et \`a condition de prendre syst\'ematiquement la polarisation orthogonale, on peut identifier structure presque CR polaris\'ee transversalement int\'egrable et feuilletage transversalement holomorphe.
\medskip
Cependant, il n'est pas possible en g\'en\'eral d'associer une structure presque CR polaris\'ee
{\it compl\`etement int\'egrable} \`a un feuilletage transversalement holomorphe $\Cal F$ donn\'e. En effet, cela suppose en plus que la r\'ealisation $E$ choisie de $N\Cal F$ soit involutive {\it dans $TX$}, ce qui n'a aucune raison d'\^etre vrai.
\medskip
Il s'agit l\`a d'un point fondamental, sur lequel nous souhaitons insister. Il suffit pour se convaincre de cette diff\'erence de comparer (3.\lastnum[-1]) avec
$$
[E^{0,1},E^{0,1}]\subset E^{0,1}
$$
qui est la condition d'int\'egrabilit\'e requise pour passer de l'int\'egrabilit\'e transverse \`a la compl\`ete int\'egrabilit\'e.
\medskip
Finissons cette section en tirant quelques cons\'equences pratiques de la proposition 3.1. Soit $T$ une structure presque CR transversalement int\'egrable. L'existence du feuilletage transversalement holomorphe $\Cal F$ signifie qu'il existe au voisinage de tout point des coordonn\'ees feuillet\'ees $(z,t)$ telles que
$$
N=T\Cal F=_{loc}\text{Vect}_\Bbb R (\partial/\partial t_1,\hdots, \partial/\partial t_d)
\leqno (3.\numerote)
$$
et 
$$
N_\Bbb C\Cal F=_{loc}\text{Vect}_\Bbb C(\partial/\partial z_1,\hdots,\partial/\partial z_p,\partial/\partial\bar z_1,\hdots,\partial/\partial \bar z_p)
\leqno (3.\numerote)
$$
formules qui utilisent implicitement les \'egalit\'es
$$
d=\dim \Cal F=\dim N\qquad\text{ et }\qquad 2p=\dim E=\text{codim }\Cal F.
\leqno (3.\numerote)
$$  
Comme la projection (3.\lastnum[-10]) est un isomorphisme entre $E_\Bbb C$ et $N_\Bbb C\Cal F$ qui, de plus, pr\'eserve (3.\lastnum[-5]) et (3.\lastnum[-9]), on en conclut qu'il existe, pour tout $i$ entre $1$ et $p$ un unique vecteur $e_i$ de $E^{0,1}$ dont la projection par $\pi$ est \'egale \`a $\partial/\partial\bar z_i$. Ecrivons plus pr\'ecis\'ement
$$
e_i=\partial/\partial\bar z_i+a_i^j \partial/\partial t_j
\leqno (3.\numerote)
$$
o\`u les $a_i^j$ sont des fonctions lisses en $(z,t)$ \`a valeurs complexes.

\remark{Notation}
Dans la formule (3.\lastnum[0]), et puis syst\'ematiquement dans la suite du texte, nous utilisons la convention d'Einstein : un indice r\'ep\'et\'e deux fois est somm\'e (et ce d'ailleurs, qu'il soit plac\'e en haut, en bas, ou altern\'e).
\endremark
\medskip
On peut en fait montrer un peu plus, lorsqu'on est en pr\'esence d'une structure CR polaris\'ee.

\proclaim{Lemme 3.3}
Soit $T$ une structure presque CR polaris\'ee transversalement int\'egra\-ble et soit $(z,t)$ un syst\`eme local de coordonn\'ees feuillet\'ees. Alors il existe un unique $p$-uplet de vecteurs $e_i$ de type {\rm (3.\lastnum[0])} formant une base locale de $E^{0,1}$.

Qui plus est, si $T$ est compl\`etement int\'egrable, les $e_i$ v\'erifient la relation
$$
[e_i,e_j]=0\qquad\text{ pour tout }i, j\text{ entre }1\text{ et } p.
\leqno (3.\numerote)
$$ 
\endproclaim

\demo{Preuve}
Il suffit de montrer la nullit\'e des crochets. Soient donc $i$ et $j$ deux entiers entre $1$ et $p$. On d\'ecompose dans la base locale $(e_k,\bar e_k,\partial/\partial t_k)$ de $T_\Bbb C X$
$$
[e_i,e_j]=\alpha_k e_k+\beta_k\bar e_k+\gamma_k\partial/\partial t_k.
\leqno (3.\numerote)
$$
Par compl\`ete int\'egrabilit\'e, le crochet $[e_i, e_j]$ appartient \`a $E^{0,1}$, donc tous les $\beta_k$ et les $\gamma_k$ sont nuls. Mais la formule (3.\lastnum[-2]) implique que ce crochet est dirig\'e le long de $N_\Bbb C$. Finalement, les $\alpha_k$ aussi doivent \^etre nuls, et (3.\lastnum[0]) se transforme en (3.\lastnum[-1]).
$\square$
\enddemo

Dans toute la suite, nous utiliserons syst\'ematiquement pour faire des calculs les coordonn\'ees feuillet\'ees et la base locale $(e_k,\bar e_k,\partial/\partial t_k)$ de $T_\Bbb C X$ associ\'ee. Par abus de notation, on \'ecrira $(z,t)\in X$ pour $M\in X$ et $M$ est repr\'esent\'e par $(z,t)$ dans un syst\`eme de coordonn\'ees feuillet\'ees.

\proclaim{Proposition 3.4}
Soit $T$ une structure CR polaris\'ee. Soit $(z,t)\in X$. Alors la forme de Levi de $E$ en $(z,t)$ est donn\'ee par la matrice
$$
\pmatrix
[\bar e_1, e_1] &\hdots &[\bar e_1, e_p]\cr
\vdots &&\vdots\cr
[\bar e_p, e_1]&\hdots &[\bar e_p, e_p]
\endpmatrix
(z,t).
\leqno (3.\numerote)
$$
\endproclaim

\demo{Preuve}
Par d\'efinition, la forme de Levi est (\`a normalisation pr\`es)
$$
(v,w)\in E^{1,0}\times E^{1,0}\longmapsto p_{N_\Bbb C}[v,\bar w]\in N_\Bbb C.
$$
Dans la base locale $\bar e_i$ de $E^{1,0}$, cela donne exactement l'expression (3.\lastnum[0]).
$\square$
\enddemo

Nous verrons dans la section suivante des exemples non Levi-plats.

\vskip.5cm
\subhead
{\bf 4. Exemples de vari\'et\'es CR polaris\'ees}
\endsubhead
\numcount=0
\example{Vari\'et\'es sasakiennes}
Rappelons qu'une vari\'et\'e lisse compacte riemannienne $(X,g)$ est sasakienne si son c\^one riemannien 
$$
\Cal C(X)=(X\times\Bbb R_{>0},r^2g+dr^2)
\leqno (4.\numerote)
$$
(o\`u $r$ est la coordonn\'ee de $\Bbb R_{>0}$) admet une structure complexe $J$ invariante par dilatations
$$
(x,r)\in\Cal C(X)\longmapsto  (x,\lambda r)\in\Cal C(X)
\leqno (4.\numerote)
$$
et k\"ahl\'erienne pour sa m\'etrique.
\medskip
En identifiant $X$ \`a l'hypersurface $X\times\{1\}$ de $\Cal C(X)$, on la munit d'une structure CR $(E, J_{\vert E})$ de codimension un. Par ailleurs, si l'on d\'efinit
$$
\xi=Jr\dfrac{\partial}{\partial r}
\leqno (4.\numerote)
$$
on v\'erifie facilement (cf. \cite{B-G} et \cite{Sp}) que $\xi$ est un champ tangent \`a $X$ dont le flot est transverse \`a $E$ et pr\'eserve la structure CR. Autrement dit, si $\Phi_t$ d\'esigne le flot de $\xi$, on a 
$$
\left (\Phi_t \right )_*E^{0,1}=E^{0,1}
\leqno (4.\numerote)
$$
et
$$
\text{Vect}_\Bbb R\ \xi\oplus E=TX.
\leqno (4.\numerote)
$$
Mais la version infinit\'esimale de (4.\lastnum[-1]) est
$$
[\xi, E^{0,1}]\subset E^{0,1}
\leqno (4.\numerote)
 $$
c'est-\`a-dire exactement l'int\'egrabilit\'e transverse de $E$ pour la polarisation
$$
N=\text{Vect}_\Bbb R\ \xi.
\leqno (4.\numerote)
$$
Combin\'ee \`a l'int\'egrabilit\'e de $E^{0,1}$, on obtient ainsi qu'une vari\'et\'e sasakienne est une vari\'et\'e CR polaris\'ee.
\medskip
Observons que cela donne de nombreux exemples non Levi-plats.
\endexample

\example{Structures CR de codimension r\'eelle une invariantes sous l'action d'un flot transverse}
Il s'agit d'une g\'en\'eralisation - sans m\'etrique k\"ahl\'erienne transverse - de l'exemple pr\'ec\'edent. On part de $(E,J)$ structure CR de codimension r\'eelle une sur $X$ et on suppose que $(E,J)$ est invariante sous l'action d'un flot transverse. Autrement dit, il existe un champ $\xi$ sur $X$ de flot $\Phi_t$ v\'erifiant (4.\lastnum[-3]) et (4.\lastnum[-2]).
\medskip
On voit alors par (4.\lastnum[-1]) que $(E,J)$ est polaris\'ee pour (4.\lastnum[0]).
\endexample

\example{Feuilletages transversalement holomorphes de codimension complexe u\-ne}
Soit $\Cal F$ un feuilletage transversalement holomorphe de codimension complexe une sur $X$. Choisissons une distribution $E$ de $TX$ isomorphe \`a $N\Cal F$. Le corollaire 3.2 entra\^{\i}ne que $E$ munie de la structure CR h\'erit\'ee de $N\Cal F$ est transversalement int\'egrable, et donc polaris\'ee par raison de dimension.
\medskip
Dans le cas particulier d'un flot transversalement holomorphe sur une vari\'et\'e de dimension $3$, on dispose d'une classification compl\`ete d'apr\`es M. Brunella et E. Ghys \cite{Br}, \cite{Gh} : fibrations de Seifert, feuilletages lin\'eaires du tore $T^3$, feuilletage stable associ\'e \`a une suspension d'un diff\'eomorphisme hyperbolique de $T^2$, suspension d'un automorphisme de $\Bbb P^1$, feuilletages transversalement affines sur $\Bbb S^2\times \Bbb S^1$ et enfin feuilletages de $\Bbb S^3$ induits par un champ de vecteurs holomorphe de $\Bbb C^2$ dans le domaine de Poincar\'e et quotients. 
\endexample

\example{Suspensions \`a base complexe}
Soit $B$ une vari\'et\'e compacte complexe. Soit $\overline B$ un rev\^etement galoisien de $B$ de groupe $\Gamma$ et soit
$$
\rho\ :\ \Gamma\longrightarrow \text{Diff}(F)
\leqno (4.\numerote)
$$
une repr\'esentation de $\Gamma$ dans les diff\'eomorphismes d'une vari\'et\'e lisse compacte $F$.
\medskip
L'action fibr\'ee
$$
(z,t)\in \overline B\times F\longmapsto (g\cdot z,(\rho(g)^{-1}(t))\in \overline B\times F
\leqno (4.\numerote)
$$
o\`u $g\in \Gamma$ agit par translations sur $\overline B$, d\'efinit par passage au quotient une vari\'et\'e
$$
X:=\overline B\times F/\langle \Gamma\rangle
\leqno (4.\numerote)
$$
qui fibre sur $B$ de fibre $F$.
\medskip
Comme $B$ est complexe, le feuilletage vertical de $\overline B\times F$ par copies de $F$ descend en un feuilletage transversalement holomorphe $\Cal F$ sur $X$. De surcro\^{\i}t, le feuilletage horizontal par copies de $\overline B$ descend en un feuilletage par vari\'et\'es complexes $\Cal G$ sur $X$. Toutes les feuilles sont des rev\^etements holomorphes de $B$. Notant $E$ le fibr\'e tangent de $\Cal G$, on voit ais\'ement que $E$ est une structure CR polaris\'ee. En fait, pour $(z,t)$ syst\`eme de coordonn\'ees locales de $\overline B\times F$, et donc de $X$, on a que
$$
e_i=\dfrac{\partial}{\partial \bar z_i}.
\leqno (4.\numerote)
$$
Ces exemples de structures CR polaris\'ees sont toutes Levi-plates.
 \endexample

\vskip.5cm
\subhead
{\bf 5. Action d'un diff\'eomorphisme proche de l'identit\'e}
\endsubhead

\numcount=0
Dans cet article, ainsi que nous l'avons signal\'e dans l'introduction et au vu des d\'efinitions de la section 2, nous travaillons toujours \`a polarisation fix\'ee. De ce fait, \'etant donn\'ee une structure presque polaris\'ee $T$, les diff\'eomorphismes $f$ que nous utiliserons en priorit\'e vont pr\'eserver la polarisation $N$ de $T$, i.e. vont v\'erifier
$$
(f_*N)_x=N_{f(x)}\qquad\text{ pour tout }x\in X.
\leqno (5.\numerote)
$$
Cependant, afin d'avoir un th\'eor\`eme d'existence d'espace local de modules dans le cadre le plus g\'en\'eral possible (cf. sections 7 et 10), nous serons \'egalement amen\'es \`a consid\'erer des diff\'eomorphismes quelconques.
\medskip
On notera Diff$_N(X)$ le groupe des diff\'eomorphismes $C^\infty$ de $X$ pr\'eservant $N$, et Diff$^l_N(X)$ le groupes des diff\'eomorphismes de classe $L^2_l$ qui compl\`ete Diff$_N(X)$ pour la norme $l$.
\medskip
Supposons $T$ compl\`etement int\'egrable. Dans ce cas, Diff$_N(X)$ est le groupe des diff\'eomorphismes $\Cal F$-feuillet\'es, c'est-\`a-dire envoyant une feuille de $\Cal F$ sur une feuille de $\Cal F$. En coordonn\'ees $(z,t)$, on a alors
$$
f(z,t)=(h(z), g(z,t)).
\leqno (5.\numerote)
$$
Soit $T'$ une d\'eformation de $T$ \`a polarisation constante cod\'ee par $\omega_\Bbb C$ (au sens de (2.9)). Soit $f$ un \'el\'ement de Diff$_N(X)$. Consid\'erons la distribution $f_*T'$ de $T_\Bbb C X$. Tenant compte de la d\'ecomposition
$$
T'=(E^{0,1})'\oplus iN
\leqno (5.\numerote)
$$
on obtient
$$
f_*T'=f_*(E^{0,1})'\oplus iN
\leqno (5.\numerote)
$$
En particulier, il ressort de (5.\lastnum[0]) que $f_* T'$ est une structure presque CR polaris\'ee par $N$. 
\medskip
Dans le cas o\`u $f$ est quelconque, la composante $h$ de (5.\lastnum[-2]) d\'epend aussi de $t$ et la distribution $f_*N$ n'est plus \'egale \`a $N$. Cependant, si $f$ est suffisamment proche de l'identit\'e, la distribution $iN$ sera encore transverse \`a $f_* (E^{0,1})'$, si bien qu'on {\it d\'efinira} $f_* T'$ par la formule (5.\lastnum[0]).
\medskip
Supposons $f$ suffisamment proche de l'identit\'e pour que $f_* T'$ soit une d\'e\-for\-ma\-ti\-on \`a polarisation constante de $T$. Codons-la par le morphisme $\alpha_\Bbb C$.
Notre objectif est de relier $\alpha_\Bbb C$ \`a $\omega_\Bbb C$.
\medskip
On se place comme d'habitude dans une carte locale feuillet\'ee et on suppose que $f$ est suffisamment proche de l'identit\'e pour \^etre \`a valeurs dans cette m\^eme carte. On peut donc \'ecrire
$$
\omega_\Bbb C=\omega_{1,0}\oplus\omega_{N_\Bbb C}=\omega_{ij}e_i^*\otimes \bar e_j\oplus n_{ij}e_i^*\otimes \partial/\partial t_j
\leqno (5.\numerote)
$$
pour $\omega_{ij}$ et $n_{ij}$ des fonctions lisses en $(z,t)$ \`a valeurs complexes et
$$
\alpha_\Bbb C=\alpha_{1,0}\oplus\alpha_{N_\Bbb C}:=\alpha_{ij}e_i^*\otimes \bar e_j\oplus p_{ij}e_i^*\otimes \partial/\partial t_j
\leqno (5.\numerote)
$$
pour $\alpha_{ij}$ et $p_{ij}$ des fonctions lisses en $(z,t)$ \`a valeurs complexes.
\medskip
Posons
$$
\partial_i=\partial/\partial z_i\quad\bar\partial_i=\partial /\partial \bar z_i\quad \partial_{t_i}=\partial/\partial t_i
\leqno (5.\numerote)
$$

\proclaim{Proposition 5.1}
Sous les hypoth\`eses et notations pr\'ec\'edentes, on a 
\medskip
\noindent (i) Si $f\in ${\rm Diff}$_N(X)$,
$$
\alpha_{lj}(\bar\partial _k \bar h_l-\omega_{ki}\partial_i\bar h_l)=\omega_{ki}\partial_i h_j-\bar\partial _k h_j 
\leqno (5.\numerote)
$$
pour $j,k=1,\hdots, p$ ; et
$$
\left\{
\eqalign{
p_{lj}(\bar\partial _k \bar h_l-\omega_{ki}\partial_i\bar h_l)&=n_{ki}\partial_{t_i}g_j-\bar\partial_k g_j-a_k^l\partial_{t_l}g_j\cr
&+\omega_{ki}(\partial_ig_j+\bar a_i^l\partial_{t_l}g_j-\partial_ih_l\bar a_l^j-\partial_i\bar h_la_l^j)\cr
&+\bar\partial_k\bar h_la_l^j+\bar\partial_k h_l\bar a_l^j
}
\right .
\leqno (5.\numerote)
$$
pour $j=1,\hdots, d$ et $k=1,\hdots, p$.
\medskip
\noindent (ii) Si $f$ est quelconque et si les coefficients $n_{ij}$ de {\rm (5.\lastnum[-4])} sont tous nuls,
$$
\left\{
\eqalign{
\alpha_{lj}(\bar\partial _k \bar h_l+a_k^m\partial_{t_m}\bar h_l
-\omega_{ki}&(\partial_i\bar h_l+\bar a_i^m\partial_{ t_m}\bar h_l))=\cr
&\omega_{ki}(\partial_i h_j+\bar a_i^m\partial_{ t_m }h_j)-(\bar\partial _k h_j+a_k^m\partial_{ t_m }h_j)
}
\right .
\leqno (5.\numerote)
$$
pour $j,k=1,\hdots, p$.
\endproclaim

\remark{Remarque}
Dans le cas $f$ quelconque, on peut bien s\^ur pousser les calculs pour donner la formule des $p_{jk}$ ainsi que les formules lorsque les $n_{ij}$ ne sont plus nuls. Nous n'avons pas jug\'e utile de les donner, car nous n'en aurons pas besoin dans cet article.
\endremark

\demo{Preuve}
(i) C'est un calcul fastidieux mais direct. De (3.13) et (5.\lastnum[-8]), on d\'eduit les formules
$$
f_*e_i=\bar\partial_i\bar h_j e_j+\bar\partial_ih_j\bar e_j+(\bar\partial_ig_j+ a_i^l\partial_{t_l}g_j-\bar\partial_i\bar h_l a_l^j-\bar\partial_ih_l\bar a_l^j)\partial/\partial t_j
\leqno (5.\numerote)
$$
et
$$
f_*\partial/\partial t_j=\partial_{t_i}g_j\partial/\partial t_j.
\leqno (5.\numerote)
$$
Soit $k$ un entier compris entre $1$ et $p$. On va calculer 
$$
f_*(e_k-\omega_\Bbb C(e_k))
\leqno (5.\numerote)
$$
de deux fa\c cons diff\'erentes. Tout d'abord, \`a partir de (5.\lastnum[-8]), on a imm\'ediatement que (5.\lastnum[0]) vaut
$$
f_*e_k-\omega_{kl}f_*\bar e_l-n_{kl}f_*\partial/\partial t_l
\leqno (5.\numerote)
$$
et on utilise ensuite (5.\lastnum[-3]) et (5.\lastnum[-2]) pour obtenir une premi\`ere 
expression explicite de (5.\lastnum[-1]), \`a savoir
$$
\eqalign{
(\bar\partial_k\bar h_j-\omega_{ki}\partial_i\bar h_j)e_j+&(\bar\partial_kh_j-\omega_{ki}\partial_ih_j)\bar e_j\cr
+&\Big (\bar\partial_kg_j+a_k^l\partial_{t_j}g_j-\bar\partial_k\bar h_la_l^j-\bar\partial_k h_l\bar a_l^j\cr
-\omega_{ki}(\partial_ig_j&+\bar a_i^l\partial_{t_l}g_j-\partial_ih_l\bar a_l^j-\partial_i\bar h_la_l^j)-n_{ki}\partial_{t_i}g_j\Big )\dfrac{\partial}{\partial t_j}.
}
\leqno (5.\numerote)
$$
Mais par ailleurs, par d\'efinition de $\alpha_\Bbb C$, et en appliquant (5.\lastnum[-11]) et (5.\lastnum[-9]), il doit exister $v=v_ie_i$ tel que (5.\lastnum[-2]) vaille
$$
v_ie_i-v_i\alpha_{il}\bar e_l-v_ip_{il}\partial/\partial t_l.
\leqno (5.\numerote)
$$
Voil\`a la deuxi\`eme expression de (5.\lastnum[-3]). En comparant (5.\lastnum[-1]) et (5.\lastnum[0]), on tombe sur (5.\lastnum[-8]) et (5.\lastnum[-7]).
\medskip
\noindent (ii) La d\'emarche est strictement la m\^eme et les calculs donnent directement (5.\lastnum[-6]). Nous omettons les d\'etails.
$\square$
\enddemo
\vfill
\eject
\head
II. D\'eformations \`a type constant.
\endhead
\vskip1cm
\numcount=0

Dans toute cette partie, on fixe une structure polaris\'ee $T$, et on utilise librement les notations associ\'ees, comme (1.2), (2.3), (2.4). 
%Dans quelques paragraphes (en particulier en section 12), $T$ est uniquement suppos\'ee %transversalement int\'egrable, mais c'est alors pr\'ecis\'e explicitement.
\vskip.5cm

\subhead
{\bf 6. Equations d'int\'egrabilit\'e}
\endsubhead

Soit $T'$ une d\'eformation \`a type constant de $T$. On la code par
$$
\omega=\omega_{ij}e_i^*\otimes \bar e_j
\leqno (6.\numerote)
$$
La proposition suivante exprime en fonction des $\omega_{ij}$ les diff\'erentes int\'egrabilit\'es de $T'$.

\proclaim{Proposition 6.1}
La structure presque CR polaris\'ee $T'$ est
\medskip
\noindent (i) transversalement int\'egrable si et seulement si, pour tout $i$ et $j$ entiers compris entre $1$ et $p$, et tout $k$ compris entre $1$ et $d$, on a
$$
((\bar\partial_j+a_j^l\partial_{t_l})\omega_{il}-(\bar\partial_i+a_i^l\partial_{t_l})\omega_{jl})\bar e_l+[\omega_{il}\bar e_l,\omega_{jl}\bar e_l]=0
\leqno (6.\numerote)
$$
et 
$$
\partial_{t_k}\omega_{ij}=0.
\leqno (6.\numerote)
$$

\noindent (ii) int\'egrable si et seulement si, pour tout $i$ et $j$ entre $1$ et $p$, on a {\rm (6.\lastnum[-1])} et
$$
\omega_{ik}[e_j,\bar e_k]=\omega_{jk}[e_i,\bar e_k].
\leqno (6.\numerote)
$$ 
\endproclaim
Ainsi $T'$ est polaris\'ee si et seulement si elle v\'erifie (6.\lastnum[-2]), (6.\lastnum[-1]) et (6.\lastnum[0]).

\demo{Preuve}
\noindent (i) Il faut v\'erifier l'involutivit\'e de $(T')^\Bbb C$. A partir de (2.4), de (3.1) et du lemme (3.3), on voit qu'il suffit de v\'erifier que
$$
L_{ij}=[e_i-\omega (e_i), e_j-\omega (e_j)]\in (T')^\Bbb C
\leqno (6.\numerote)
$$
et 
$$
M_{ik}=[e_i-\omega (e_i),\partial/\partial t_k]\in N_\Bbb C.
\leqno (6.\numerote)
$$
Maintenant (6.\lastnum[-1]) est v\'erifi\'ee si et seulement si 
$$
L_{ij}^{1,0}=-\omega(L_{ij}^{0,1})
$$
ce qui, en utilisant (3.14), donne 
$$
L_{ij}^{1,0}=(e_j\cdot (\omega (e_i))-e_i\cdot (\omega (e_j))+[\omega (e_i),\omega (e_j)]=0
\leqno (6.\numerote)
$$
Il suffit maintenant de d\'evelopper (6.\lastnum[0]) en tenant compte de (3.13) et (6.\lastnum[-6]) pour obtenir (6.\lastnum[-5]).
\medskip
De m\^eme, (6.\lastnum[-2]) est \'equivalent \`a 
$$
M_{ik}^{1,0}=M_{ik}^{0,1}=0
$$
soit, apr\`es d\'eveloppement, exactement (6.\lastnum[-4]).
\medskip
\noindent (ii) On veut maintenant que
$$
L_{ij}=[e_i-\omega (e_i), e_j-\omega (e_j)]\in (E')^{0,1}
\leqno (6.\numerote)
$$
ce qui entra\^{\i}ne (6.\lastnum[-1]), et donc (6.\lastnum[-6]) ; mais aussi
$$
L_{ij}^{N_\Bbb C}=0.
\leqno (6.\numerote)
$$
Un calcul direct donne alors
$$
\eqalign{
L_{ij}^{N_\Bbb C}&=[e_j,\omega (e_i)]^{N_\Bbb C}-[e_i,\omega (e_j)]^{N_\Bbb C}\cr
&=\omega_{ik}(e_j\cdot \bar e_k)-\omega (e_i)\cdot e_j-\omega_{jk}(e_i\cdot \bar e_k)+\omega (e_j)\cdot e_i\cr
&=\omega_{ik}[e_j,\bar e_k]-\omega_{jk}[e_i,\bar e_k]
}
\leqno (6.\numerote)
$$
et finalement (6.\lastnum[-6]).
$\square$
\enddemo
\vskip.5cm
\subhead
{\bf 7. Formes isochrones}
\endsubhead

\numcount=0
Compte-tenu de (7.3), il est tentant de travailler avec des formes $\omega$ ind\'ependantes de $t$. Le lemme suivant montre que c'est possible.

\proclaim{Lemme 7.1}
Soient $(z,t)$ des coordonn\'ees locales d\'efinies sur $U\subset X$. Soit $\omega$ une $1$-forme complexe sur $E^{0,1}$ \`a valeurs dans $E^{1,0}$. On suppose que les coefficients de $\omega$ dans la base $e_i^*\otimes\bar e_j$ (cf. {\rm (7.1)}) sont ind\'ependants de $t$.

Alors, dans tout autre syst\`eme de coordonn\'ees locales feuillet\'ees $(w,s)$ sur $U$, les coefficients de $\omega$ dans la base $f_i^*\otimes\bar f_j$ associ\'ee sont ind\'ependants de $s$.
\endproclaim

\demo{Preuve}
C'est un calcul sans difficult\'e. Les deux syst\`emes de coordonn\'ees sont reli\'es par
$$
(w,s)=(h(z), g(z,t))
\leqno (7.\numerote)
$$
avec $h$ holomorphe et $g$ de classe $C^\infty$. Partant de l'\'ecriture (3.13) des $e_i$, et appliquant la jacobienne du changement de cartes (7.\lastnum[0]), on obtient
$$
e_i=\bar\partial_i\bar h_j\partial/\partial\bar w_j+(\bar\partial_ig_j+a_i^k\partial g_j/\partial t_k)\partial/\partial s_j.
\leqno (7.\numerote)
$$
Mais l'unicit\'e du lemme 3.3 entra\^{\i}ne que l'on a
$$
e_i=\bar\partial_i \bar h_j f_j
\leqno (7.\numerote)
$$
d'o\`u, par dualit\'e,
$$
e_i^*=\bar\partial _j \bar h_i^{-1} f_j^*
\leqno (7.\numerote)
$$
et finalement
$$
\omega=\omega_{ij}e_i^*\otimes \bar e_j=\omega_{ij}\bar\partial_k\bar h_i^{-1}\partial_jh_l f_k ^*\otimes \bar f_l.
\leqno (7.\numerote)
$$
Comme $h$ ne d\'epend pas de $t$, ni $h^{-1}$ de $s$, on voit que les coefficients de $\omega$ dans la base $e_i^*\otimes\bar e_j$ sont ind\'ependants de $t$ si et seulement si ses coefficients dans la base $f_k^*\otimes f_l$ sont ind\'ependants de $s$.
$\square$
\enddemo
\remark{Remarque}
Insistons sur le fait que ce sont les {\it coefficients} de $\omega$ dans la base $e_i^*\otimes \bar e_j$ qui sont ind\'ependants de $t$, pas les $e_i^*\otimes \bar e_j$ eux-m\^emes, qui d\'ependent en g\'en\'eral de $t$ via les coefficients $a_i^j$ de (3.13).
\endremark
\medskip 
Plus g\'en\'eralement, si $\omega$ est cette fois une $q$-forme sur $E^{0,1}$ \`a valeurs dans $E^{1,0}$, on a localement
$$
\omega=\omega_{i_1\hdots i_q,j}e_{i_1}^*\wedge \hdots\wedge e_{i_q}^*\otimes \bar e_j
\leqno (7.\numerote)
$$
et dans un autre syst\`eme de coordonn\'ees (7.\lastnum[-5]),
$$
\omega=\omega_{i_1\hdots i_q,j}\bar\partial_{k_1}\bar h_{i_1}^{-1}\hdots \bar\partial_{k_q}\bar h_{i_q}^{-1}\partial_j h_l f_{k_1}^*\wedge\hdots\wedge f_{k_q}^*\otimes\bar f_l.
\leqno (7.\numerote)
$$
De m\^eme que dans la preuve du lemme 7.1, comme $h$ ne d\'epend pas de $t$, ni $h^{-1}$ de $s$, on voit que les coefficients de $\omega$ dans la base $e$ sont ind\'ependants de $t$ si et seulement si ses coefficients dans la base $f$ sont ind\'ependants de $s$.

\definition{D\'efinition}
Soit $\omega$ une $q$-forme sur $E^{0,1}$ \`a valeurs dans $E^{1,0}$. On dira que $\omega$ est {\it isochrone} si, dans toute \'ecriture locale {\rm (7.\lastnum[-1])}, les coefficients $\omega_{i_1\hdots i_q,j}$ sont ind\'e\-pen\-dants de $t$.
\enddefinition

On notera $A^q$ le $\Bbb C$-espace vectoriel des $q$-formes isochrones lisses sur $E^{0,1}$ \`a valeurs dans $E^{1,0}$ et, comme d'habitude, $A^q_l$ son compl\'et\'e pour la norme $l$.
\medskip
Ainsi l'espace $A^1$ repr\'esente l'espace des d\'eformations isochrones de $T$ \`a type constant, c'est-\`a-dire qui v\'erifient l'\'equation (6.3). L'espace $A^0$, l'espace des champs isochrones de $E^{1,0}$, va quant \`a lui jouer un r\^ole particulier, puisqu'il va nous permettre de coder les diff\'eomorphismes dont nous allons nous servir. 
\medskip
Dans la preuve classique de Kuranishi de construction d'un espace versel de d\'eformations pour les vari\'et\'es compactes complexes, on utilise le fait que le groupe Diff$^l(X)$ est une vari\'et\'e de Banach model\'ee sur l'espace des champs de vecteurs de classe $L^2_l$ de $X$. De plus, si $\exp$ d\'esigne l'exponentielle d'une connexion lin\'eaire fix\'ee sur $X$, l'application
$$
e\ :\ \xi\in W_l\subset \Omega^0_l(T)\longmapsto \exp (\Re\xi)\in\text{Diff}^l(X)
\leqno (7.\numerote)
$$
est une carte locale explicite en l'identit\'e pour $W_l$ voisinage suffisamment petit de $0$.
\medskip
Dans notre cas, le plus naturel serait de travailler avec les diff\'eomorphismes fixant $E$ et $N$. Mais c'est beaucoup trop restrictif en g\'en\'eral et la relation d'\'equivalence induite sur les structures CR est beaucoup trop fine pour esp\'erer un espace de modules de dimension finie.
\medskip
Le deuxi\`eme choix naturel est de travailler avec les diff\'eomorphismes pr\'eservant la polarisation $N$. Mais pour cela, nous avons besoin d'une application (7.\lastnum[0]) qui {\it en plus} v\'erifie 
$$
e(A^0_l\cap W_l)\subset\text{Diff}_N^l(X).
\leqno (7.\numerote)
$$
Nous ne savons pas si une telle application existe en g\'en\'eral. Ceci est reli\'e au probl\`eme de savoir si Diff$_N^l(X)$ est une vari\'et\'e Banachique.
\medskip
Nous allons au contraire utiliser une connexion li\-n\'e\-ai\-re quel\-con\-que sur $X$. L'ima\-ge de 
$A^0_l\cap W_l$ par $e$ va d\'efinir un sous-ensemble de diff\'eomorphismes qui ne respectent ni $E$ ni $N$ au sens strict, mais les pr\'eservent \`a l'ordre 1. On voit en effet que les champs de $A^0$ sont tangents \`a $E$ et pr\'eservent $N$. Une fois qu'on applique $e$, le d\'efaut de pr\'eservation de ces distributions provient donc uniquement de la "non-invariance transverse" de la connexion lin\'eaire. Toutefois, nous verrons au chapitre 12 que, pour $\chi\in A^0_l$ et $s\in\Bbb R$, le terme d'ordre 1 en $s$ de $e(s\chi)$ pr\'eserve $E$ et $N$.
\medskip
Soit $\omega$ une d\'eformation \`a type constant de $T$. \`A un diff\'eomorphisme $f$, on associe la d\'eformation \`a type constant
$$
f\cdot \omega:=P_{iso}(f_*\omega)^{1,0}
\leqno (7.\numerote)
$$
i.e. on fait agir $f$ sur $\omega$ comme dans la proposition 5.1 (notre $f_*\omega$ n'est rien d'autre que le $\alpha_\Bbb C$ de la proposition 5.1) ; mais nous ne gardons que la partie $(1,0)$ de $f_*\omega$, puis nous projetons orthogonalement sur le sous-espace des formes isochrones afin de retomber sur une d\'eformation \`a type constant isochrone.
\medskip
Le probl\`eme de cette d\'efinition est qu'elle ne pr\'eserve ni la condition d'int\'egrabi\-li\-t\'e (6.2) ni la condition (6.4). Nous verrons en section 10 comment y rem\'edier.
\medskip
Avant d'aller plus loin, consid\'erons un cas particulier classique important, o\`u nous pourrons donner une interpr\'etation plus g\'eom\'etrique des choses.

\definition{D\'efinition}
Nous dirons que $E$ poss\`ede une {\it connexion invariante par holonomie (CIH)} si $E$ poss\`ede une connexion lin\'eaire qui induit sur $N\Cal F$, via l'isomorphisme (3.2), une connexion lin\'eaire invariante par le pseudo-groupe d'holonomie de $\Cal F$.
\enddefinition

Dans ces conditions, consid\'erons une connexion obtenue sur $TX$ en sommant une CIH sur $E$ et une connexion lin\'eaire quelconque sur $T\Cal F$. L'exponentielle associ\'ee \`a une telle connexion permet de d\'efinir une application $e$ qui v\'erifie (7.\lastnum[-1]).

\remark{Remarque}
Le premier et le dernier exemple de la section 4 poss\`edent une CIH. Dans le deuxi\`eme, c'est aussi le cas si les orbites du flot sont ferm\'ees.
\endremark
\medskip

Il est toutefois important de noter que, contrairement au cas classique de \cite{Ku2}, les diff\'eomorphismes obtenus ainsi ne recouvrent pas un voisinage de l'identit\'e dans Diff$_N(X)$. En effet, dans l'\'ecriture locale
$$
e(\chi)(z,t)= (h(z),g(z,t))
\leqno (7.\numerote)
$$
la partie $g$ ne d\'epend pas de $\chi$ mais uniquement des coefficients $a_i^j$.
\medskip
Sous cette hypoth\`ese, la formule (7.\lastnum[-1]) se simplifie. On a en effet

\proclaim{Lemme 7.2}
Soit $f\in\text{Diff}_N(X)$ et soit $\omega\in A^1$. Alors la d\'eformation \`a type constant
$(f_*\omega ) ^{1,0}$ est isochrone.
\endproclaim

\demo{Preuve}
La forme $(f_*\omega ) ^{1,0}$ est donn\'ee localement par la formule (5.8). Comme $\omega$ est isochrone et que $h$ (une des composantes de $f$, cf. (5.2)) est ind\'ependante de $t$, on voit que les coefficients de $f\cdot\omega$ sont ind\'ependants de $t$, i.e. que cette forme est isochrone.
$\square$
\enddemo

Ainsi la projection $P_{iso}$ dans la formule (7.\lastnum[-1]) est inutile.
G\'eom\'etriquement, cela signifie la chose suivante. Nous ne regardons pas l'action de $f$ sur $E$, puisqu'elle ne pr\'eserve pas en g\'en\'eral $E$ (il faudrait sinon travailler avec le groupe des diff\'e\-o\-mor\-phis\-mes pr\'eservant $N$ et $E$ ce qui poserait bien des probl\`emes dans la suite - en particulier car il est souvent trop petit) ; nous regardons l'action de $f$ sur le fibr\'e quotient $TX/N$, identifi\'e \`a $E$ via la projection naturelle (3.2). En d'autres termes, m\^eme si $f$ ne pr\'eserve pas $E$, il pr\'eserve $N\Cal F$ et modifie la structure complexe sur $N\Cal F$. Mais comme nous avons une identification fixe entre $E$ et $N\Cal F$, tout cela permet de d\'efinir sans ambigu\"{\i}t\'e la d\'eformation \`a type constant induite par $f$.
\medskip
Revenons maintenant au cadre g\'en\'eral.
En combinant (7.\lastnum[-1]) et (7.\lastnum[-3]), on obtient une application $\Bbb C$-analytique
$$
(\chi,\omega)\in A^0_{l+1}\cap W_{l+1}\times A^1_l\longmapsto e(\chi)\cdot\omega=P_{iso}(e(\chi)_*\omega)^{1,0}\in A^1_l.
\leqno (7.\numerote)
$$
On remarque en effet que l'application $e$ est $\Bbb C$-analytique de m\^eme que l'application $*$ (cf. \cite{Ku2}), la projection sur $A^1_l$ et celle sur les formes isochrones.
\vskip.5cm
\subhead
{\bf 8. L'op\'erateur $\bar\partial$}
\endsubhead
\numcount=0

Soit $T'$ une d\'eformation \`a type constant de $T$ cod\'ee par $\omega$. Supposons que $T'$ soit transversalement int\'egrable. Alors $\omega$ est isochrone. De plus l'\'equation (6.2) se simplifie en 
$$
(\bar\partial_j\omega_{il}-\bar\partial_i\omega_{jl})\bar e_l+[\omega_{il}\bar e_l,\omega_{jl}\bar e_l]=0
\leqno (8.\numerote)
$$
qu'il est tentant, par comparaison avec le cas classique, de simplifier en
$$
\bar\partial \omega+[\omega,\omega]=0.
\leqno (8.\numerote)
$$
Ceci signifie d'abord que l'on pose
$$
[\omega, \alpha]:=-\dfrac{1}{2}e_i^*\wedge e_j^*\otimes [\omega (e_i),\alpha (e_j)]
\leqno (8.\numerote)
$$
Ainsi (8.\lastnum[0]) appliqu\'ee \`a la paire de vecteurs $(e_j,e_i)$ donne (8.\lastnum[-1]).

\remark{Remarque}
Ce crochet est \'egal \`a $-1/2$-fois celui de Kuranishi, d'o\`u la diff\'erence dans les formules d'int\'egrabilit\'e.
\endremark
\medskip

Ceci suppose ensuite qu'on ait en coordonn\'ees locales (en utilisant (6.1))
$$
\bar\partial\omega :=\bar\partial_j\omega_{il}e_j^*\wedge e_i^*\otimes \bar e_l.
\leqno (8.\numerote)
$$
On peut \'etendre facilement cette formule \`a tout 
$$
\chi=\chi_i\bar e_i\in A^0
\leqno (8.\numerote)
$$ 
en \'ecrivant en coordonn\'ees locales
$$
\bar \partial\chi:=\bar\partial_j\chi_i e_j^*\otimes \bar e_i
\leqno (8.\numerote)
$$
et \`a tout $\omega\in A^q$, d\'ecompos\'e selon (8.6)
$$
\bar\partial\omega:=\bar\partial_j\omega_{i_1\hdots i_q,l}e_j^*\wedge e_{i_1}^*\wedge\hdots\wedge e_{i_q}^*\otimes \bar e_l.
\leqno (8.\numerote)
$$
On a
\proclaim{Lemme 8.1}
Les \'equations {\rm (8.\lastnum[-3]), (8.\lastnum [-1])} et {\rm (8.\lastnum[0])} d\'efinissent globalement des o\-p\'e\-ra\-teurs 
$$
\bar\partial\ :\ A^q\longrightarrow A^{q+1}.
\leqno (8.\numerote)
$$

De plus, ces op\'erateurs forment une suite exacte
$$
0\aro >>>H^0 \aro >>> A^0\aro >\bar\partial >>A^1\aro >\bar\partial>>\hdots\aro >\bar\partial >> A^p\aro >>> 0
\leqno (8.\numerote)
$$
o\`u $H^0$ est l'espace vectoriel des champs holomorphes de $E^{1,0}$.
\endproclaim 

\demo{Preuve}
V\'erifions par calcul que $\bar\partial$ est globalement d\'efini pour les champs. Pour cela, on se place dans un nouveau syst\`eme de coordonn\'ees (7.1). Soit $\chi$ un \'el\'ement de $A^0$. On tire de (8.\lastnum[-4]) et (7.3) la formule
$$
\chi=\chi_i\partial_i h_j\bar f_j
\leqno (8.\numerote)
$$
d'o\`u, par (8.\lastnum[-4])
$$
\eqalign{
\bar\partial\chi=&\bar\partial_l(\chi_i\partial_ih_j)f_l^*\otimes \bar f_j\cr
=&\bar\partial_m\chi_i\bar\partial_l\bar h_m^{-1}\partial_ih_jf_l^*\otimes \bar f_j
}
\leqno (8.\numerote)
$$
en utilisant le fait que $h$ est holomorphe. Mais (7.5) implique que cette formule co\"{\i}ncide avec (8.\lastnum[-5]). Ceci montre que $\bar\partial$ est globalement d\'efini de $A^0$ dans $A^1$.
\medskip
Pour \'etendre la d\'efinition de $\bar\partial$ \`a $A^q$ pour $q>0$, il suffit d'utiliser la formule de Cartan. Ainsi, pour une 1-forme $\omega$, on peut d\'efinir globalement $\bar\partial$ en posant, pour $\xi$ et $\eta$ des champs de $E^{0,1}$,
$$
\bar\partial\omega (\xi,\eta)=-(\bar\partial (\omega(\xi)))\eta+(\bar\partial (\omega(\eta)))\xi-\omega([\xi,\eta])
\leqno (8.\numerote)
$$
Mais on v\'erifie imm\'ediatement que cette formule, appliqu\'e localement \`a $(e_i,e_j)$ donne exactement (8.\lastnum[-8]) (\`a cause de (3.14)). Autrement dit, les d\'efinitions locale (8.\lastnum[-8]) et globale (8.\lastnum[0]) sont identiques, donc le $\bar\partial$ de (8.\lastnum[-8]) est bien un op\'erateur global de $A^1$ dans $A^2$. 
\medskip
Par induction, en comparant la formule de Cartan au rang $q$ \`a la formule (8.\lastnum[-5]), on montre de m\^eme que (8.\lastnum[-5]) est bien un op\'erateur global de $A^q$ dans $A^{q+1}$.
\medskip
Enfin, puisque $\bar\partial$ v\'erifie la formule de Cartan, on a automatiquement
$$
\bar\partial\circ \bar\partial \equiv 0
\leqno (8.\numerote)
$$
et on en d\'eduit ais\'ement la suite exacte (8.\lastnum[-4]).
$\square$
\enddemo

Observons que les op\'erateurs (8.\lastnum[-5]), de par leur d\'efinition locale, s'\'etendent en des op\'erateurs continus
$$
\bar\partial\ :\ A^q_l\longrightarrow A^{q+1}_{l-1}.
\leqno (8.\numerote)
$$
Rappelons maintenant que les espaces $A^q_l$ sont des espaces de Hilbert. On peut donc associer \`a $\bar\partial$ son adjoint de Hilbert
$$
\bar\partial^*\ :\ A^{q+1}_{l-1}\longrightarrow A^q_l
\leqno (8.\numerote)
$$
d\'efini par
$$
(\bar\partial \omega,\alpha)_{l-1}=(\omega,\bar\partial^*\alpha)_l
\leqno (8.\numerote)
$$
pour tout $\omega\in A^q_l$ et $\alpha\in A^{q+1}_{l-1}$.

\remark{Remarque}
Nous voulons insister sur le fait que nous utilisons ici, et dans toute la suite, l'adjoint de Hilbert, et non pas l'adjoint des op\'erateurs diff\'erentiels. D'ail\-leurs, tel que nous l'avons d\'efini, $\bar\partial$ n'est m\^eme pas un op\'erateur diff\'erentiel, puisqu'il est d\'efini sur $A^q$ qui n'est pas l'espace des sections lisses d'un fibr\'e vectoriel sur $X$, mais l'espace des sections lisses {\it isochrones} d'un tel fibr\'e, cf. la discussion en d\'ebut de section 9.
\endremark
\vskip.5cm
\subhead
{\bf 9. L'op\'erateur $\bar D$}
\endsubhead

\numcount=0
Les espaces $A^q$ d\'efinis en section 8 ne sont pas des espaces de sections d'un fibr\'e vectoriel. En fait, soit
$$
\Omega^q(E^{0,1})\otimes E^{1,0}
\leqno (9.\numerote)
$$
le fibr\'e des $q$-formes sur $E^{0,1}$ \`a valeurs dans $E^{1,0}$. On voit que $A^q$ est un sous-espace strict de l'espace des sections de $\Omega^q(E^{0,1})\otimes E^{1,0}$. Plus pr\'ecis\'ement, il s'agit du sous-espace des sections isochrones. 
\medskip
De ce fait, l'op\'erateur $\bar\partial$ de la section 8, d\'efini sur les $A^q$ n'est pas un op\'erateur diff\'erentiel. On peut toutefois rem\'edier \`a ce probl\`eme en l'\'etendant en un op\'erateur diff\'erentiel agissant sur les sections de $\Omega^q(E^{0,1})\otimes E^{1,0}$.
Pour cela, on d\'enote par $A^q(t)$ ces espaces de sections et on pose pour
$$
\chi=\chi_i\bar e_i\in A^0(t)
\leqno (9.\numerote)
$$
l'op\'erateur
$$
\bar\partial \chi=(e_j\cdot \chi_i)e_j^*\otimes \bar e_i\in A^1(t)
\leqno (9.\numerote)
$$
qu'on \'etend par la formule de Cartan (8.12) aux $p$-formes. Observons que (9.\lastnum[0]) redonne (8.6) lorsque $\chi$ est isochrone.

\proclaim{Lemme 9.1}
L'op\'erateur d\'efini en {\rm (9.\lastnum[0])} est globalement d\'efini et induit une suite exacte
$$
0\aro >>>H^0(t) \aro >>> A^0(t)\aro >\bar\partial >>A^1(t)\aro >\bar\partial>>\hdots\aro >\bar\partial >> A^p(t)\aro >>> 0
\leqno (9.\numerote)
$$
o\`u $H^0(t)$ est l'espace vectoriel des champs $C^\infty$ de $E^{1,0}$ holomorphes en la variable $z$.
\endproclaim

\demo{Preuve}
C'est une adaptation directe du lemme 8.1.
$\square$
\enddemo

L'op\'erateur ainsi \'etendu est alors un op\'erateur diff\'erentiel d'ordre $1$. Toutefois, il n'est pas elliptique, et le complexe (9.\lastnum[0]) n'est pas non plus elliptique.
\medskip
On peut pourtant associer \`a $\bar\partial$ un op\'erateur induisant un complexe elliptique. Pour cela, on consid\`ere le fibr\'e $\Omega^q(T^\Bbb C)\otimes E^{1,0}$ des $q$-formes sur $T^{\Bbb C}$ \`a valeurs dans $E^{1,0}$. Soit $\Cal A^q$ l'espace des sections globales lisses de 
$\Omega^q(T^\Bbb C)\otimes E^{1,0}$. On d\'efinit localement, pour 
$$
\chi=\chi_i\bar e_i\in \Cal A^0=A^0(t)
\leqno (9.\numerote)
$$
l'op\'erateur
$$
\bar D\chi:=\bar\partial_k\chi_ie_k^*\otimes \bar e_i+\partial_{t_k}\chi_idt_k\otimes\bar e_i.
\leqno (9.\numerote)
$$
Le premier terme du membre de droite n'est rien d'autre que le $\bar\partial$ d\'efini en (9.\lastnum[-3]).

\proclaim{Lemme 9.2}
L'op\'erateur local $\Bar D$ de la formule {\rm (9.\lastnum[0])} est un op\'erateur globalement d\'efini.
\endproclaim

\demo{Preuve}
Elle est tout-\`a-fait similaire \`a celle du lemme 8.1. On se place dans un nouveau syst\`eme de coordonn\'ees (8.1). Soit $\chi$ un \'el\'ement de $\Cal A^0$ d\'ecompos\'e comme en (9.\lastnum[-1]). En utilisant le fait que $\bar f_j^*$ est nulle sur $T^\Bbb C$, on trouve que
$$
dt_k=\bar\partial_jg^{inv}_kf_j^*+\partial_{s_j} g^{inv}_k ds_j
\leqno (9.\numerote)
$$
o\`u
$$
(z,t)=(h^{-1}(w), g^{inv}(w,s)).
\leqno (9.\numerote)
$$
Coupl\'e \`a (8.10), (9.\lastnum[-1]) permet de calculer
$$
\eqalign{
\bar D(\chi)=&\bar\partial_k(\chi_i\partial_ih_j)f_k^*\otimes\bar f_j+\partial_{s_k}(\chi_i\partial_ih_j)ds_k^*\otimes \bar f_j\cr
=&\bar\partial_m\chi_i(\bar\partial_kh_m^{-1}f_k^*\otimes \partial_ih_j\bar f_j)\cr
&+\partial_{t_m} \chi_i (\bar\partial_kg_m^{inv}f_k^*+\partial_{s_k}g_m^{inv}ds_k^*)\otimes\partial_ih_j\bar f_j
}
\leqno (9.\numerote)
$$
ce qui donne, en tenant compte de (7.3), (7.4) et (9.\lastnum[-2]), exactement la d\'efinition (9.\lastnum[-3]).
$\square$
\enddemo

Comme dans la preuve du lemme 8.1, on \'etend ensuite $\bar D$ aux autres $\Cal A^q$ en utilisant la formule de Cartan et on obtient une suite exacte
$$
0\aro >>> H^0\aro >>>\Cal A^0\aro >\bar D >>\Cal A^1\aro >\bar D>>\hdots\aro >\bar D >> \Cal A^p\aro >>> 0.
\leqno (9.\numerote)
$$

\proclaim{Lemme 9.3}
Le complexe {\rm (9.\lastnum[0])} est un complexe elliptique.
\endproclaim

\demo{Preuve}
On calcule la suite associ\'ee de symboles. Soit $(z,t,v)\in T^* X$. En utilisant la d\'ecomposition
$$
v=v^{1,0}\oplus v^{0,1}\oplus v^{N_\Bbb C}\in E^{1,0}\oplus E^{0,1}\oplus N_\Bbb C=T_\Bbb CX
\leqno (9.\numerote)
$$
on v\'erifie que le symbole de $\bar D$ appliqu\'e \`a un \'el\'ement $\omega$ de $\Omega^q(T^\Bbb C)\otimes E^{1,0}$ en $(z,t,v)$ vaut
$$
(v^{1,0}+v^{N_\Bbb C})\wedge \omega.
\leqno (9.\numerote)
$$
On v\'erifie facilement, \`a partir de (9.\lastnum[0]), l'exactitude de la suite (9.\lastnum[-2]).
$\square$
\enddemo

On notera que l'op\'erateur $\bar D$ restreint \`a $A^q$ est exactement l'op\'erateur $\bar\partial$.

\remark{Remarque}
On notera par contre que l'op\'erateur local $\partial_t=\bar D-\bar \partial$ n'est pas bien d\'efini globalement. Ceci est li\'e au fait qu'on n'a pas en g\'en\'eral de d\'ecomposition locale $d=\partial+\bar\partial+\partial_t$.
\endremark

\vskip.5cm
\subhead
{\bf 10. Mod\`ele local de l'espace des d\'eformations \`a type constant trans\-ver\-sa\-lement int\'egrables et de l'espace des feuilletages transversalement holomorphes \`a type diff\'erentiable fix\'e}
\endsubhead

\numcount=0

Nous sommes maintenant en mesure de d\'ecrire localement l'espace des d\'efor\-ma\-tions \`a type constant transversalement int\'egrables et celui des feuilletages transversalement holomorphes \`a type diff\'erentiable fix\'e.
\medskip
Soit $\Cal T_{tc}^l$ l'ensemble des structures presque CR polaris\'ees transversalement in\-t\'e\-gra\-bles de classe $L^2_l$ \`a $E$ et \`a $N$ fix\'es. C'est un sous-espace $\Bbb C$-analytique de la vari\'et\'e de Banach $\Cal E^l_{tc}$. Via le corollaire 2.2 et la proposition 6.1, il s'identifie localement en $T$ \`a un voisinage de $0$ dans 
$$
\{\omega\in A^1_l\quad\vert\quad \bar\partial \omega+[\omega,\omega]=0\}.
\leqno (10.\numerote)
$$
D\'ecomposons
$$
A^0_l=H^0\oplus (A^0_l)^\perp .
\leqno (10.\numerote)
$$
On pose
$$
\left\{
\eqalign{
\Phi_l\ :\ (\chi, \omega)&\in (A^0_{l+1})^\perp\cap W_{l+1}\times A^1_l\longmapsto\cr & (\bar\partial (e(\chi)\cdot\omega) +[(e(\chi)\cdot\omega,e(\chi)\cdot\omega],\bar\partial ^*\omega)\in A^2_{l-1}\times A^0_{l+1}}
\right .
\leqno (10.\numerote)
$$
et 
$$
K_l=P_1\circ\Phi_l^{-1} (0,0)
\leqno (10.\numerote)
$$
o\`u $P_1$ est la projection naturelle de $(A^0_{l+1})^\perp\times A^1_l$ sur le facteur $A^1_l$.
Posons enfin
$$
\Xi_l\ :\ (\chi, \omega)\in (A^0_{l+1})^\perp\cap W_{l+1}\times K_l\longmapsto e(\chi)\cdot \omega\in \Cal T^l_{tc}.
\leqno (10.\numerote)
$$

Nous pouvons \'enoncer le th\'eor\`eme de structure locale de $\Cal T^l_{tc}$.

\proclaim{Th\'eor\`eme 10.1}
On se place sous les hypoth\`eses et les notations pr\'ec\'edentes. Alors,

\noindent (i) Au voisinage de $0$, l'espace $K_l$ est un ensemble analytique de dimension finie.

\noindent (ii) L'application $\Xi_l$ est un isomorphisme $\Bbb C$-analytique local en $(0,0)$.
\endproclaim

Ainsi un voisinage de $0$ dans $\Cal T_{tc}^l$, c'est-\`a-dire un voisinage de $T$ dans l'espace analytique des d\'eformations \`a type constant de $T$ transversalement int\'egrables, est isomorphe analytiquement \`a un voisinage de $(0,0)$ dans le produit de $K_l$ par l'espace vectoriel  $(A^0_l)^\perp$. 
\medskip
On remarquera que nous ne pouvons pas donner des \'equations plus pr\'ecises pour $K_l$ en g\'en\'eral. En effet, nous sommes oblig\'es d'utiliser la pr\'esentation indirecte (10.\lastnum[-1]) car l'action (7.10) ne pr\'eserve pas l'int\'egrabilit\'e transverse. Cela ne change rien sur le plan {\it th\'eorique} mais rend beaucoup plus d\'elicat le calcul d'exemples.
\medskip
Dans le cas particulier o\`u $T$ poss\`ede une CIH, on peut \^etre plus pr\'ecis.

\proclaim{Corollaire 10.2}
Supposons de plus que $T$ poss\`ede une CIH. Alors on a 
$$
K_l=\{\omega\in A^1_l\quad\vert\quad \bar\partial^*\omega=\bar\partial\omega+[\omega,\omega]=0\}.
\leqno (10.\numerote)
$$ 
\endproclaim

De surcro\^{\i}t, rappelons que Diff$^{l+1}(X)$ agit analytiquement sur $\Cal T^l_{tc}$ par (7.10) et que, si $\omega$ est un point de $K_l$, l'image $\Xi_l( (A^0_{l+1})^\perp\cap W_{l+1},\omega)$ est incluse dans l'orbite de $\omega$ pour cette action. On a donc

\proclaim{Corollaire 10.3}
L'espace analytique $K_l$ est une section locale analytique en $T$ de l'action de {\rm Diff}$^{l+1}(X)$ sur $\Cal T^l_{tc}$
\endproclaim

Par section locale analytique en $T$, nous entendons un espace analytique plong\'e qui rencontre
{\it toutes les orbites de l'action passant suffisamment pr\`es de $T$ en au moins un point} (mais pas forc\'ement unique).
\medskip
Ainsi toute structure CR transversalement int\'egrable proche de $0$ est isomorphe (au sens de (7.10)) \`a un point de $K_l$. Plus encore,

\proclaim{Corollaire 10.4}
Soit $L$ un espace analytique et soit $(\omega_t)_{t\in L}$ une famille analytique de d\'eformations $L^2_l$ transversalement int\'egrables \`a type constant de $T$, i.e. l'application
$$
t\in L\longmapsto \omega_t\in A^1_l
\leqno (10.\numerote)
$$
est analytique. 

Alors il existe un germe d'application analytique en $0\in L$ (point tel que $\omega_0=0$)
$$
f\ :\ (L,0)\longrightarrow (K_l,0)
\leqno (10.\numerote)
$$
telle que les germes de $L$ en $0$ d'une part et du pull-back $f^*K_l$ en $0$ d'autre part soient isomorphes.
\endproclaim

\demo{Preuve}
C'est une cons\'equence imm\'ediate du th\'eor\`eme 10.1 et du fait que $\Xi_l$ est analytique.
$\square$
\enddemo

Il faut voir ce corollaire comme un \'equivalent de la propri\'et\'e classique de com\-pl\'e\-tude de la th\'eorie des d\'eformations de Kodaira-Spencer. Toutefois, il y a ici une diff\'erence fondamentale. Dans la th\'eorie de Kodaira-Spencer, on {\it d\'efinit} une d\'eformation comme un morphisme plat entre espaces analytiques. La base est l'espace de param\`etres et les fibres repr\'esentent la m\^eme vari\'et\'e lisse munie de structures complexes diff\'erentes. On montre {\it ensuite} - et cela occupe, par exemple, une bonne partie du livre \cite{Ku3} - que cette d\'efinition est \'equivalente \`a une d\'efinition de d\'eformation type (10.\lastnum[-1]). Autrement dit, pour que le corollaire 10.4 soit vraiment un r\'esultat de compl\'etude au sens de Kodaira-Spencer, il faudrait d\'evelopper une notion de d\'eformations type morphisme transplat (cf. \cite{Me2}) et en faire une \'etude approfondie.
\medskip
Une telle \'etude est peut \^etre possible, mais nous remarquons que personne n'a pu la r\'ealiser bien que quarante ans soient pass\'es depuis les travaux de Kodaira-Spencer... {\it Notre point de vue est justement qu'on peut obtenir un r\'esultat de compl\'etude sans avoir \`a d\'evelopper toute la machinerie technique inh\'erente \`a l'ap\-pro\-che de Kodaira-Spencer}.
\medskip
La preuve du th\'eor\`eme 10.1 est donn\'ee en section 12, \`a l'aide des pr\'eliminaires de la section 11.
\medskip
 Enfin, notons que le contenu du th\'eor\`eme 10.1 peut \^etre transpos\'e au cas des feuilletages transversalement holomorphes. En effet, d'apr\`es le corollaire 3.2, il existe une bijection entre $\Cal T_{tc}$ et l'ensemble des feuilletages transversalement holomorphes tangents \`a $N$. Convenons de dire qu'un tel feuilletage est de classe $L^2_l$ si la structure transverse est de classe $L^2_l$ (le feuilletage \'etant quant \`a lui suppos\'e $C^\infty$). Ainsi l'ensemble des feuilletages $L^2_l$ s'identifie \`a $\Cal T^l_{tc}$.
\medskip
On peut maintenant \'enoncer

\proclaim{Th\'eor\`eme 10.5}

\noindent (i) Soit $\Cal F$ un feuilletage transversalement holomorphe de classe $L^2_l$ sur $X$. On suppose qu'il existe une structure CR polaris\'ee $T$ sur $X$ de feuilletage associ\'e $\Cal F$. Alors l'ensemble des feuilletages transversalement holomorphes de classe $L^2_l$ sur $X$ proches de $\Cal F$ et diff\'eomorphes \`a $\Cal F$ forme un $\Bbb C$-espace analytique localement isomorphe \`a $K_l\times (A^0_{l+1})^\perp$.

\noindent (ii) Si, de plus, $T$ admet une CIH, alors deux points ayant m\^eme projection sur $K_l$ sont isomorphes en tant que feuilletages transversalement holomorphes.
\endproclaim

\demo{Preuve}
Le corollaire 3.2 affirme que  l'ensemble des feuilletages transversalement holomorphes de classe $L^2_l$ sur $X$ diff\'eomorphes \`a $\Cal F$ est en bijection avec $\Cal T^2_{tc}$. On applique le th\'eor\`eme 10.1. Enfin la partie (ii) est simplement une reformulation des consid\'erations du paragraphe 7.
$\square$
\enddemo
Qui plus est, on retrouve la propri\'et\'e de compl\'etude 10.4. Ainsi $K_l$ est un espace "complet" pour le feuilletage transversalement holomorphe $\Cal F$.
\medskip
Il est int\'eressant de comparer ce th\'eor\`eme avec les r\'esultats classiques de \cite{EK-N} et \cite{Gi} sur les d\'eformations de feuilletages transversalement holomorphes \`a type diff\'erentiable fixe. Dans \cite{Gi} - qui contient les \'enonc\'es les plus g\'en\'eraux -, on montre un th\'eor\`eme de versalit\'e sous les hypoth\`eses suivantes :
\medskip
\noindent i) Il existe une CIH sur le fibr\'e normal.

\noindent ii) Un certain op\'erateur diff\'erentiel est \`a image ferm\'ee.
\medskip
On voit donc que notre th\'eor\`eme peut \^etre consid\'er\'e comme un analogue du r\'esultat de versalit\'e cit\'e. D'une part, nous pouvons nous passer des hypoth\`eses i) et ii), mais au prix de l'utilisation d'une relation d'\'equivalence sur les feuilletages qui n'est pas la plus naturelle. Et m\^eme dans le cas plus g\'eom\'etrique o\`u il existe une CIH, nous pouvons nous passer de l'hypoth\`ese ii), ce qui repr\'esente un gain ind\'eniable, ce type d'hypoth\`ese \'etant tr\`es difficile \`a v\'erifier en pratique. Mais d'autre part, nous n'avons pas la versalit\'e, simplement la compl\'etude, et nous devons supposer que le feuilletage base admette une structure CR polaris\'ee.
\medskip
En ce sens, les deux \'enonc\'es ne sont pas exactement comparables. Les preuves, bien que suivant le m\^eme sch\'ema g\'en\'eral, sont d'ailleurs assez diff\'erentes et ne peuvent s'obtenir l'une de l'autre.
\medskip
L'absence de versalit\'e n'est en outre pas sans cons\'equence. Il en r\'esulte que nous ne savons pas si le germe de $K_l$ est {\it unique}. Ainsi la d\'ependance de $K_l$ en $l$ n'est pas claire (voir cependant le corollaire 15.3), de m\^eme que la d\'ependance en $T$ de l'\'enonc\'e 10.5.
\medskip
Pour r\'esumer, il faut donc comprendre les th\'eor\`emes 10.1 et 10.5 comme des \'enonc\'es affirmant l'existence d'un espace complet de {\it dimension finie}, c'est-\`a-dire comme des \'enonc\'es de {\it rigidit\'e}.

\vskip .5cm
\subhead
{\bf 11. L'op\'erateur $\Delta$}
\endsubhead
\numcount =0

Classiquement, les constructions d'espaces de modules locaux "\`a la Kuranishi" n\'ecessite l'introduction d'un op\'erateur elliptique dont le noyau est l'espace tangent \`a l'espace en question, sa dimension finie provenant ainsi directement de l'ellipticit\'e de l'op\'erateur.
\medskip
Dans notre cas, cela supposerait d'utiliser le laplacien associ\'e \`a $\bar D$. Cependant, ce choix entra\^{\i}ne beaucoup de probl\`emes. En effet, ce laplacien ne va pas en g\'en\'eral respecter les espaces $A^q$. Cela provient du fait qu'il faut pour le construire choisir une m\'etrique riemannienne sur $X$ ; et qu'il faudrait pour qu'il pr\'eserve les formes isochrones, que la m\'etrique soit elle-m\^eme isochrone, c'est-\`a-dire invariante par l'holonomie du feuilletage transversalement holomorphe. Mais une telle m\'etrique n'existe pas forc\'ement, il faut supposer le feuilletage riemannien. C'est l'hypoth\`ese faite par exemple dans \cite{EK-N} o\`u le m\^eme probl\`eme se pose (cf. section 10).
\medskip
Pour se passer de cette hypoth\`ese suppl\'ementaire, nous allons travailler avec l'op\'erateur
$$
\Delta:=\bar\partial \bar\partial^*+\bar\partial^*\bar\partial\ :\ A^q_l\longrightarrow A^q_l
\leqno (11.\numerote)
$$
qui, rappelons-le, est construit \`a partir de l'adjoint {\it hilbertien} de $\bar\partial$. Ce n'est pas un op\'erateur diff\'erentiel, juste un op\'erateur lin\'eaire entre espaces de Hilbert. Nous allons montrer qu'il est Fredholm, notion bien connue pour \^etre en quelque sorte l'\'equivalent pour les op\'erateurs lin\'eaires de l'ellipticit\'e pour les op\'erateurs diff\'erentiels. Bien s\^ur, cela va nous priver des r\'esultats de r\'egularit\'e des op\'erateurs diff\'erentiels elliptiques, ce qui n'est pas sans cons\'equence pour passer des structures de classe $L^2_l$ aux structures $C^\infty$, mais globalement cela ne va changer grand chose.

\proclaim{Proposition 11.1}
L'op\'erateur $\Delta$ est un op\'erateur de Fredholm auto-adjoint.
\endproclaim

\demo{Preuve}
L'op\'erateur $\Delta$ est de fa\c con \'evidente auto-adjoint, si bien qu'il suffit de montrer que son noyau est de dimension finie. Ceci implique en effet que son conoyau est de dimension finie, et donc que son image est ferm\'ee (cf. \cite{Pa}).
\medskip
Maintenant, $\omega\in A^q_l$ est dans le noyau de $\Delta$ si et seulement si
$$
\bar\partial\omega=\bar\partial^*\omega=0.
\leqno (11.\numerote)
$$
Soit $\Omega^q_{inv}$ le faisceau des germes de $q$-formes sur $E^{0,1}$ \`a valeurs dans $E^{1,0}$, qui soient isochrones et $\bar\partial$-ferm\'ees. L'\'equation (11.\lastnum[0]) entra\^{\i}ne
$$
\text{Ker }\Delta\subset H^0(X,\Omega^q_{inv}).
\leqno (11.\numerote)
$$
Il suffit donc de montrer que ce groupe est de dimension finie. Ce qui est le contenu du lemme 11.2 suivant, qui conclut donc la preuve.
$\square$
\enddemo

\proclaim{Lemme 11.2}
Les groupes de cohomologie $H^0(X,\Omega^q_{inv})$ sont de dimension finie.
\endproclaim

\demo{Preuve}
C'est un r\'esultat classique, qu'on peut trouver dans \cite{D-K} ou \cite{G-M}. Alternativement, on le retrouve facilement en notant que, si $\Omega^q$ d\'esigne le faisceau des germes de $q$-formes sur $T^\Bbb C$ \`a valeurs dans $E^{1,0}$, qui sont isochrones et $\bar\partial$-ferm\'ees, alors l'inclusion naturelle
$$
\Omega^q_{inv}\subset \Omega^q
\leqno (11.\numerote)
$$
induit une inclusion des groupes de cohomologie
$$
H^0(X,\Omega^q_{inv})\subset H^0(X,\Omega^q).
\leqno (11.\numerote)
$$
Or, on obtient, \`a partir de (9.10), une r\'esolution
$$
0\aro >>> \Omega^q\aro >>>\Omega^q(T^\Bbb C)\otimes E^{1,0}\aro >\bar D>>
\Omega^{q+1}(T^\Bbb C)\otimes E^{1,0}\aro >\bar D>>\hdots
\leqno (11.\numerote)
$$
qui est elliptique par le lemme 9.2.
$\square$
\enddemo

\remark{Remarque}
Il est bien connu que les groupes $H^i(X,\Omega^q_{inv})$ peuvent \^etre de dimension infinie pour $i>0$.
\endremark
\medskip
Comme cons\'equence de la proposition 11.1, il existe un op\'erateur 
$$
G\ :\ A^q_l=\text{Ker } \Delta\oplus \text{Im }\Delta\longrightarrow \text{Im }\Delta
\leqno (11.\numerote)
$$
tel que l'on ait
$$
G\Delta(\omega)=\Delta G(\omega)=\omega_1\quad\text{ pour }\omega=\omega_0\oplus\omega_1\in \text{Ker }\Delta\oplus \text{Im }\Delta.
\leqno (11.\numerote)
$$
\vskip.5cm
\subhead
{\bf 12. Preuve du th\'eor\`eme 10.1}
\endsubhead

\numcount=0

{\it Nous prouvons d'abord le (ii)}. L'action de Diff$^{l+1}_N(X)$ sur la vari\'et\'e de Banach $A^1_l$ est analytique (cf. \cite{Do1}). De plus l'application $e$ d\'efinie en (7.8) - essentiellement l'exponentielle pour une m\'etrique riemannienne sur $X$ - est analytique. Il en est donc de m\^eme de $\Xi_l$. Etendons $\Xi_l$ \`a
$$
\tilde\Xi_l\ :\ (\chi, \omega)\in (A^0_{l+1})^\perp\cap W_{l+1}\times \text{Ker }\bar\partial^*\longmapsto e(\chi)\cdot \omega\in A^1_{l}.
\leqno (12.\numerote)
$$
 Calculons sa diff\'erentielle en $(0,0)$. On a

\proclaim{Lemme 12.1}
La diff\'erentielle de $\tilde\Xi_l$ en $(0,0)$ est \'egale \`a 
$$
d_{(0,0)}\tilde\Xi_l(\xi,\omega)=\omega-\bar\partial \xi
\leqno (12.\numerote)
$$
\endproclaim

\demo{Preuve}
On se place en cartes locales. Sachant que, d'apr\`es (3.13) et (7.8), on a, pour $s$ petit,
$$
e(s\chi_ie_i)(z,t)\simeq (z_1+s\chi_1,\hdots, z_p+s\chi_p,t_1+s\chi_i\bar a_{i1},\hdots, t_d+s\chi_i\bar a_{id})
\leqno (12.\numerote)
$$
o\`u le signe $\simeq$ signifie \'egalit\'e \`a l'ordre un en $s$, on en d\'eduit que la formule (5.10), sous l'hypoth\`ese
$$
\alpha_s:=(e(s\chi)_* s\omega )^{1,0}
\leqno (12.\numerote)
$$
s'\'ecrit en cartes locales \`a l'ordre $1$ en $s$
$$
(e(s\chi)_* s\omega )^{1,0}\simeq s\omega-s\bar\partial\chi \ .
\leqno (12.\numerote)
$$
Comme $P_{iso}$ est lin\'eaire, elle pr\'eserve la partie d'ordre un en $s$ et on a
$$
e(s\chi)\cdot s\omega \simeq P_{iso}(s\omega-s\bar\partial\chi )\ .
\leqno (12.\numerote)
$$
Mais $\omega$ et $\chi$ sont isochrones donc $P_{iso}$ agit ici comme l'identit\'e. En prenant $s=1$, on obtient donc bien que 
la diff\'erentielle de $\Xi_l$ en $(0,0)$ est (12.\lastnum[-4]).
$\square$
\enddemo

Il suffit maintenant de montrer que (12.\lastnum[-4]) est inversible. Il s'ensuivra que $\tilde\Xi_l$ est un isomorphisme analytique local, d'apr\`es la version du th\'eor\`eme d'inversion locale de \cite{Do2}. Et par restriction aux structures transversalement int\'egrables, il en sera de m\^eme de $\Xi_l$. La preuve du th\'eor\`eme 10.1, (ii) se termine donc avec le

\proclaim{Lemme 12.2}
La diff\'erentielle {\rm(12.\lastnum[-4])} est inversible.
\endproclaim

\demo{Preuve}
Consid\'erons l'application
$$
\Gamma\ :\ \alpha\in A^1_{l}\longmapsto (-G \bar\partial^*\alpha,\alpha-\bar\partial G \bar\partial^*\alpha)\in A^0_{l+1}\times A^1_l.
\leqno (12.\numerote)
$$
On voit que la premi\`ere composante appartient en fait \`a Im $\Delta$, c'est-\`a-dire \`a $(A^0_{l+1})^\perp$. De plus, on a
$$
\eqalign{
\bar\partial^*(\alpha-\bar\partial G \bar\partial^*\alpha)=&\bar\partial^*\alpha-\Delta G\bar\partial^*\alpha\cr
=&\Delta G\bar\partial^*\alpha- \Delta G\bar\partial^*\alpha=0
}
\leqno (12.\numerote)
$$
car Im $\bar\partial^*$ est orthogonal au noyau de $\Delta$ (qui n'est autre que le noyau de $\bar\partial$ puisque nous sommes sur les $0$-formes), donc inclus dans l'image de $\Delta$, ce qui conclut par (11.8).
\medskip
D\`es lors, la deuxi\`eme composante de $\Gamma$ est \`a valeurs dans Ker $\bar\partial^*$. Enfin, on a
$$
\eqalign{
\Gamma\circ d_{(0,0)}\tilde\Xi_l(\xi,\alpha)=&\Gamma (\alpha-\bar\partial \xi)\cr
=&(G\Delta\xi, \alpha-\bar\partial\xi+\bar\partial G\Delta\xi)\cr
=&(\xi,\alpha)
}
\leqno (12.\numerote)
$$
 ce qui ach\`eve la preuve.
$\square$
\enddemo  

{\it Nous prouvons maintenant le (i)}. Il ne reste plus qu'\`a montrer que $K_l$ est de dimension finie. Il r\'esulte du lemme 12.1 que la diff\'erentielle de l'application $\Phi_l$ en $(0,0)$ est \'egale \`a
$$
(\chi,\omega)\in (A^0_{l+1})^\perp\times A^1_l\longmapsto (-\bar\partial\omega,\bar\partial^*\omega)\in A^2_{l-1}\times A^0_{l+1}.
\leqno (12.\numerote)
$$
On observe que cette diff\'erentielle ne d\'epend pas de $\chi$. D\`es lors, en projetant par $P_1$ on obtient que l'espace tangent de Zariski de $K_l$ en $0$ est le noyau de l'op\'erateur $\Delta$ appliqu\'e aux $1$-formes, qui est de dimension finie par la proposition 11.1.

\vskip .5cm
\subhead
{\bf 13. Mod\`ele local de l'espace des d\'eformations \`a type constant com\-pl\`e\-te\-ment int\'egrables}
\endsubhead
\numcount=0

Soit $\Cal C_{tc}^l$ l'ensemble des structures CR polaris\'ees de classe $L^2_l$ \`a $E$ et \`a $N$ fix\'es. C'est un sous-espace $\Bbb C$-analytique de la vari\'et\'e de Banach $\Cal E_{tc}$. Via le corollaire 2.2 et la proposition 6.1, il s'identifie localement en $T$ \`a un voisinage de $0$ dans 
$$
\{\omega\in A^1_l\quad\vert\quad \bar\partial \omega+[\omega,\omega]=0\text{ et v\'erifie (6.4)}\}.
\leqno (13.\numerote)
$$
Remarquons que, si l'on pose
$$
F_l=\{\omega\in A^1_l\quad\vert\quad \omega\text{ satisfait (6.4)}\}
\leqno (13.\numerote)
$$
on d\'efinit ainsi un sous-espace vectoriel ferm\'e de $A^1_l$ tel que, via les identifications (13.\lastnum[-1]) et (10.1),
$$
\Cal C_{tc}^l=\Cal T_{tc}^l\cap F_l.
\leqno (13.\numerote)
$$
Posons enfin
$$
K^0_l:=\pi_{K_l}\circ\Xi_l^{-1}(F_l)=\pi_{K_l}\circ\Xi_l^{-1} (\Cal C_{tc}^l)
\leqno (13.\numerote)
$$
o\`u $\pi_{K_l}$ d\'esigne la projection sur $K_l$ ; et
$$
\Xi^0_l:=(\Xi_l)\vert _{((A^0_l)^\perp\times K^0_l)}.
\leqno (13.\numerote)
$$

Le th\'eor\`eme de structure locale de $\Cal T_{tc}^l$ s'obtient alors comme cons\'equence directe du th\'eor\`eme 10.1.

\proclaim{Th\'eor\`eme 13.1}

\noindent (i) Au voisinage de $0$, l'espace $K^0_l$ est un ensemble analytique de dimension finie.

\noindent (ii) L'application $\Xi^0_l$ est un isomorphisme $\Bbb C$-analytique local en $(0,0)$.
\endproclaim

\demo{Preuve}
Par (10.5), l'espace $K^0_l$ est la projection sur le facteur $K_l$ de l'espace analytique $\Xi_l^{-1}(F_l)$. Il s'agit donc d'un sous-espace analytique de $K_l$.
$\square$
\enddemo

\remark{Remarque}
Il n'est pas vrai que $K^0_l$ soit donn\'e comme l'intersection de $K_l$ et de $F_l$. Ceci vient de ce que l'action des diff\'eomorphismes que nous avons d\'efinie ne pr\'eserve pas l'int\'egrabilit\'e compl\`ete, mais simplement l'int\'egrabilit\'e transverse. Une fois encore, ceci ne change pas grand chose au niveau {\it th\'eorique}, mais pose des probl\`emes pratiques : cela rend tr\`es difficile le {\it calcul explicite} de $K^0_l$ sur des exemples.
\endremark
\medskip
On dispose aussi des deux m\^emes corollaires qu'en section 10.

\proclaim{Corollaire 13.2}
L'espace analytique $K_l^0$ est une section locale analytique en $T$ de l'action de Diff$^{l+1}_N(X)$ sur $\Cal C^l_{tc}$
\endproclaim

Ainsi toute structure CR polaris\'ee proche de $0$ est isomorphe \`a un point de $K_l$. Plus encore,

\proclaim{Corollaire 13.3}
Soit $L$ un espace analytique et soit $(\omega_t)_{t\in L}$ une famille analytique de d\'eformations $L^2_l$ compl\`etement int\'egrables \`a type constant de $T$, i.e. l'application
$$
t\in L\longmapsto \omega_t\in A^1_l
\leqno (13.\numerote)
$$
est analytique. 

Alors il existe un germe d'application analytique en $0\in L$ (point tel que $\omega_0=0$)
$$
f\ :\ (L,0)\longrightarrow (K^0_l,0)
\leqno (13.\numerote)
$$
telle que les germes de $L$ en $0$ d'une part et du pull-back $f^*K^0_l$ en $0$ d'autre part soient isomorphes.
\endproclaim

\demo{Preuve}
C'est une cons\'equence imm\'ediate du th\'eor\`eme 13.1 et du fait que $\Xi^0_l$ est analytique.
$\square$
\enddemo

L\`a encore, il faut voir ce corollaire comme un \'equivalent de la propri\'et\'e classique de compl\'etude de la th\'eorie des d\'eformations de Kodaira-Spencer, avec les r\'eserves expliqu\'ees au paragraphe 10.

\vskip .5cm
\subhead
{\bf 14. Le cas Levi-plat}
\endsubhead
\numcount=0

Dans cette section, nous supposons que $T$ est Levi-plate. On peut alors simplifier les r\'esultats pr\'ec\'edents. En effet, d'une part, l'\'equation d'int\'egrabilit\'e (6.4) est automatiquement v\'erifi\'ee par application de la proposition 3.4. On a donc

\proclaim{Lemme 14.1} 
Soit $T'$ une d\'eformation \`a type constant de $T$.
Alors $T'$ est transversalement int\'egrable si et seulement si $T'$ est compl\`etement int\'egrable.
\endproclaim

D'autre part, comme cons\'equence directe de ce lemme, on voit que les espaces $K_l$ et $K_l^0$ sont identiques.

\proclaim{Proposition 14.2}
Supposons $T$ Levi-plate. Alors les mod\`eles locaux $K_l$ des d\'eformations \`a type constant transversalement int\'egrables et $K_l^0$ des d\'eformations \`a type constant polaris\'ees sont identiques.
\endproclaim

Ainsi le th\'eor\`eme 13.1 se r\'eduit au th\'eor\`eme 10.1 dans le cas Levi-plat. Si l'on regarde de plus pr\`es la d\'efinition des diff\'eomorphismes au paragraphe 7, on voit que, si l'on suppose que la connexion lin\'eaire utilis\'ee est somme d'une connexion sur $E$ et d'une connexion sur $N$, alors les diff\'eomorphismes de $e(A^0_l)$ appartiennent \`a la vari\'et\'e banachique Diff$_E^{0,l}(X)$ des diff\'eomorphismes de classe $L_2^l$ qui pr\'eservent $E$ et {\it fixent chaque feuille de $E$}. D\`es lors, les formules (7.10) et (7.12) se simplifient en
$$
e(\chi)\cdot\omega=P_{iso}(e(\chi)_*\omega).
\leqno (14.\numerote)
$$

\medskip
Allons plus loin. La pr\'esence de deux feuilletages transverses entra\^{\i}ne l'existence d'un atlas produit.

\proclaim{Proposition 14.3}
Supposons $T$ Levi-plate. Soit $\Cal G$ le feuilletage Levi-plat. Alors on peut trouver un atlas feuillet\'e pour $\Cal F$ et $\Cal G$ simultan\'ement, i.e. un atlas dont les changements de cartes sont du type
$$
(w,s)=(h(z), g(t))
\leqno (14.\numerote)
$$
\endproclaim

\demo{Preuve} 
On d\'efinit $\Cal F$ par une collection de submersions
$$
f_i\ :\ U_i\longrightarrow \Bbb C^p
\leqno (14.\numerote)
$$
telles qu'il existe des biholomorphismes
$$
\phi_{ij}\ :\ f_i(U_i\cap U_j)\subset\Bbb C^p\longrightarrow f_j(U_i\cap U_j)\subset\Bbb C^p
\leqno (14.\numerote)
$$
v\'erifiant
$$
f_j\equiv \phi_{ij}\circ f_i\quad\text{ sur } U_i\cap U_j.
\leqno (14.\numerote)
$$
De m\^eme, d\'efinissons $\Cal G$ par une collection de submersions
$$
g_i \ :\ U_i\longrightarrow \Bbb R^d
\leqno (14.\numerote)
$$
v\'erifiant
$$
g_j\equiv \psi_{ij}\circ g_i\quad\text{ sur } U_i\cap U_j.
\leqno (14.\numerote)
$$
pour des diff\'eomorphismes $\psi_{ij}$ bien choisis.
\medskip
Il suffit maintenant de prendre comme cartes feuillet\'ees 
$$
x\in U_i\longmapsto (f_i(x),g_i(x))=(z,t)\in\Bbb C^p\times\Bbb R^d
\leqno (14.\numerote)
$$
et on voit que les changements de cartes seront du type
$$
(w,s)=(\phi_{ij}(z),\psi_{ij}(t))
\leqno (14.\numerote)
$$
donc du type (14.\lastnum[-7]).
$\square$
\enddemo

\definition{D\'efinition}
On appellera un atlas v\'erifiant {\rm (14.\lastnum[-7])} un atlas {\it isochrone}. Les cartes qu'il contient seront aussi dites {\it isochrones}.
\enddefinition

Maintenant ceci implique imm\'ediatement le

\proclaim{Corollaire 14.4}
Dans une carte feuillet\'ee isochrone, on a 
$$
e_i=\dfrac{\partial}{\partial \bar z_i}.
\leqno (14.\numerote)
$$
\endproclaim

Supposons maintenant que $T$ poss\`ede une CIH. Alors, pour $\chi\in A^0$ suffisamment petit, on d\'eduit de (14.\lastnum[0]) que, dans une carte feuillet\'ee isochrone,
$$
e(\chi)(z,t)=(h(z), t).
\leqno (14.\numerote)
$$
Ceci a pour cons\'equence que (7.9) devient
$$
e(A^0_l\cap W_l)\subset \text{Diff}_{E,N}^{0,l}(X)
\leqno (14.\numerote)
$$
o\`u $\text{Diff}_{E,N}^{0,l}$ d\'esigne le groupe des diff\'eomorphismes de $X$ de classe $L_2^l$ qui pr\'eservent $E$ et $N$ et {\it fixent les feuilles de $E$} (c'est le sens de l'exposant $0$ dans la notation). Donc (7.10) est en fait
$$
e(\chi)\cdot \omega=e(\chi)_*\omega
\leqno (14.\numerote)
$$
pour $\omega$ dans $A^1$.
\medskip
En r\'esum\'e,

\proclaim{Corollaire 14.5}
Supposons $T$ Levi-plate et admettant une CIH. Alors 
\medskip
\noindent (i) L'ap\-pli\-ca\-tion $e$ est \`a valeurs dans les diff\'eomorphismes de $X$ fixant $E$ et $N$ et chaque feuille de $E$.

\noindent (ii) L'action {\rm (7.10)} co\"{\i}ncide avec l'action standard des diff\'eomorphismes.
\endproclaim

Ainsi, dans le cas Levi-plat avec CIH, on identifie deux structures polaris\'ees proches si elles sont CR-isomorphes et non pas si elles induisent le m\^eme feuilletage transversalement holomorphe.

\vskip.5cm
\subhead
{\bf 15. Le cas des suspensions \`a base complexe}
\endsubhead
\numcount=0

Dans cette situation, on va pouvoir donner une tr\`es jolie description de l'espace $K_l$. Commen\c cons par observer le

\proclaim{Lemme 15.1}
Soit $X$ une suspension \`a base complexe $B$. Alors la structure polaris\'ee induite admet une CIH.
\endproclaim

\demo{Preuve}
Notons $\pi$ la submersion de $X$ sur $B$. Il d\'ecoule de la construction rappel\'ee en section 4 que 
$$
E\simeq\pi^* TB.
\leqno (15.\numerote)
$$
Il suffit maintenant de prendre une connexion quelconque sur $TB$ et de la relever \`a $E$ via (15.\lastnum[0]).
$\square$
\enddemo

Enon\c cons le

\proclaim{Th\'eor\`eme 15.2}
Soit $X$ une submersion \`a base complexe. Soit $K_l$ l'espace local de structures polaris\'ees proches de $T$. 

Alors $K_l$ s'identifie naturellement \`a l'espace de Kuranishi de la base $B$.
\endproclaim

Le mot "naturellement" prendra tout son sens dans la preuve. Mais avant de passer \`a celle-ci, donnons un exemple frappant d'application du th\'eor\`eme.

\example{Feuilletages lin\'eaires sur le tore $\boldkey T^3$}
Soit $\Bbb E_\tau$ la courbe elliptique de module $\tau\in\Bbb H$. Soit $\Gamma$ le groupe d'automorphismes du rev\^etement
$$
\Bbb C^*\longrightarrow \Bbb E_\tau.
\leqno (15.\numerote)
$$
Fixons un r\'eel $\alpha$ et d\'efinissons la repr\'esentation 
$$
\rho \ :\ p\in\Bbb Z\simeq\Gamma\longmapsto \bigg ([x]\in\Bbb R/\Bbb Z\mapsto [x+p\alpha]\in\Bbb R/\Bbb Z\bigg )\in \text{Diff}(\Bbb S^1).
\leqno (15.\numerote)
$$
La suspension 
$$
X=\Bbb E_\tau\times_\rho \Bbb S^1
\leqno (15.\numerote)
$$
est un tore r\'eel $T^3$ muni d'un feuilletage $\Cal F$ par cercles et d'un feuilletage transverse $\Cal G$ dit lin\'eaire. Si $\alpha$ est un rationnel $p/q$, alors les feuilles de $\Cal G$ sont toutes compactes et isomorphes \`a la m\^eme courbe elliptique $\Bbb E_\sigma$, un rev\^etement de degr\'e $q$ de $\Bbb E_\tau$. Si $\alpha$ est irrationnel, toutes les feuilles de $\Cal G$ sont denses et isomorphes \`a $\Bbb C^*$.
\medskip
Par application du th\'eor\`eme 15.2, l'espace $K_l$ de d\'eformations \`a type constant s'identifie \`a l'espace de Kuranishi de $\Bbb E_\tau$, c'est-\`a-dire \`a un voisinage de $\tau$ dans le demi-plan de Poincar\'e $\Bbb H$.
\medskip
Le point remarquable est que {\it cet espace de d\'eformations est totalement ind\'e\-pen\-dant de la pente $\alpha$}.
\medskip
Il est instructif de comparer cet espace avec celui des d\'eformations de la structure CR $E$ \`a type diff\'erentiable fix\'e, i.e. de l'ensemble des structures CR Levi-plates proches de $E$ et induisant le m\^eme feuilletage. 
\medskip
Lorsque $\alpha$ est rationnel, alors cet espace s'identifie \`a l'espace des lacets $C^\infty$ \`a valeurs dans $\Bbb H$ proches du lacet constant \'egal \`a $\sigma$ (cf. \cite{Me2}).
\medskip
 On trouve dans \cite{Sl} des calculs de cohomologie qui montrent que, dans le cas irrationnel, cet espace d\'epend des propri\'et\'es arithm\'etiques de $\alpha$. Il est de dimension $1$ lorsque $\alpha$ v\'erifie une condition diophantienne. Lorsque cette condition n'est pas v\'erifi\'ee, on ne conna\^{\i}t pas cet espace, mais les calculs de \cite{Sl} effectu\'es dans un cas tr\`es proche sugg\`erent fortement qu'il est de dimension infinie. 
\medskip
Pour comprendre le lien entre \cite{Sl} et le th\'eor\`eme 15.2, il suffit de se rappeler que nous ne consid\'erons que les structures CR {\it polaris\'ees} proches de $E$ et induisant le m\^eme feuilletage. Ainsi l'espace des structures polaris\'ees est un sous-espace parfois strict de l'espace des structures CR Levi-plates proches \`a type diff\'erentiable fix\'e.
\medskip
Plus pr\'ecis\'ement, tenant compte  de la proposition 3.1, on voit que, dans le cas rationnel, le fait d'\^etre polaris\'ee implique que toutes les feuilles de $\Cal G$ doivent \^etre isomorphes \`a la m\^eme courbe elliptique. Dans l'espace des lacets pr\'ec\'edemment d\'ecrit, on ne voit donc que les lacets constants proches de $\sigma$. Plus encore, le fait d'\^etre polaris\'e entra\^{\i}ne que la feuille doit admettre un quotient fini d'ordre $q$ (le d\'enominateur de la pente $\alpha$). En d'autres termes, on ne d\'eforme pas la feuille $\Bbb E_\sigma$ arbitrairement, mais on d\'eforme le rev\^etement $\Bbb E_\sigma\to\Bbb E_\tau$ de degr\'e $q$. Voil\`a pourquoi l'on tombe \`a la fin sur un voisinage de $\tau$ dans $\Bbb H$.
\medskip
Enfin, dans le cas irrationnel, les structures polaris\'ees sont \'egalement les structures CR qui sont invariantes par le groupe $\Gamma$, donc descendent \`a la courbe $\Bbb E_\tau$.
\endexample

\demo{Preuve du th\'eor\`eme 15.2}
Appelons comme d'habitude $J$ l'op\'erateur CR et appelons $J_B$ l'op\'erateur complexe de $B$. La submersion $\pi$ est une submersion CR au sens o\`u elle envoie chaque feuille de $\Cal G$ holomorphiquement sur $B$. On peut donc choisir un mod\`ele diff\'erentiable dans lequel 
$$
J=\left (\pi_{\vert \Cal G}\right )^*J_B
\leqno (15.\numerote)
$$
ce qu'on peut reformuler en 
$$
\pi^*E_B^{0,1}=E^{0,1}
\leqno (15.\numerote)
$$
o\`u $E_B^{0,1}$ est le sous-fibr\'e des champs $(0,1)$ de $J_B$.
\medskip
On note $A^p_B$ l'espace des $p$-formes sur $E_B^{0,1}$ \`a valeurs dans $E_B^{1,0}$. Il r\'esulte de (15.\lastnum[0]) que
$$
A^p=\pi^*A^p_B.
\leqno (15.\numerote)
$$
En effet, une $p$-forme isochrone sur $X$ se rel\`eve en une $p$-forme isochrone $\Gamma$-invariante sur $\overline B\times F$. Mais sur un tel produit, le caract\`ere isochrone signifie la constance sur les fibres $F$. Et une $p$-forme $\Gamma$-invariante et constante sur les fibres $F$ n'est rien d'autre que le relev\'e d'une $p$-forme sur $B$.
\medskip
Par ailleurs l'op\'erateur $\bar\partial$ descend lui aussi comme op\'erateur $\bar\partial$ de $B$ ; et son dual fait de m\^eme.
\medskip
On voit donc au final que
$$
\pi_*(K_l)=\{\omega\in A^1_{B,l}\quad\vert\quad\bar\partial\omega+[\omega,\omega]=\bar\partial^*\omega=0\}.
\leqno (15.\numerote) 
$$
On serait alors tent\'es de conclure directement qu'il s'agit de l'espace de Kuranishi de $B$ en comparant (15.\lastnum[0]) avec les d\'efinitions formellement identiques de \cite{Ku2} et \cite{Ku3}.
\medskip
Il y a toutefois un petit probl\`eme. Dans la d\'efinition de l'espace de Kuranishi l'op\'erateur $\bar\partial^*$ est le dual de l'op\'erateur $\bar\partial$ {\it en tant qu'op\'erateur diff\'erentiel}. Or dans notre cas, nous avons utilis\'e (cf. (8.15), (8.16) et la remarque qui suit) l'adjoint {\it de Hilbert}. Ce ne sont pas du tout les m\^emes op\'erateurs!
\medskip
Pour surmonter ce probl\`eme, notons tout d'abord que la formule (14.10) montre que les diff\'eomorphismes employ\'es se projettent sur Diff$^l (B)$. En tenant compte des corollaires 13.3 et 14.5, on en d\'eduit que l'espace $\pi_*K_l$ est un espace de d\'eformations de $B$ complet en $J_B$.
\medskip
Mais alors la th\'eorie classique de Kodaira-Spencer permet de conclure que $\pi_*K_l$ est isomorphe \`a l'espace de Kuranishi de $B$ si nous arrivons \`a montrer que la dimension de son espace tangent de Zariski en $J_B$ est \'egale \`a la dimension du premier groupe de cohomologie $H^1(B,\Theta)$ de $B$ \`a valeurs dans le faisceau $\Theta$ des champs holomorphes.
\medskip
Ce groupe peut se calculer, via les isomorphismes de Dolbeault et en tenant compte de (15.\lastnum[-1]), en prenant le quotient du noyau de l'op\'erateur $\bar\partial$ sur $A^1$ par l'image de $\bar\partial$ sur $A^0$ ; ou encore, ce qui ne change rien, par  le quotient du noyau de l'op\'erateur $\bar\partial$ sur $A^1_l$ par l'image de $\bar\partial$ sur $A^0_{l+1}$.
\medskip
Maintenant, on montre facilement, gr\^ace \`a (11.8) et (11.2), que le noyau de l'op\'erateur $\Delta$ sur $A^1_l$ co\"{\i}ncide exactement avec ce quotient. Or d'apr\`es (12.10), le noyau de $\Delta$ est la dimension de Zariski de $K_l$. Ce qui termine la preuve.
$\square$
\enddemo

\proclaim{Corollaire 15.3}
Sous les hypoth\`eses du th\'eor\`eme 15.2, l'espace $K_l$ est ind\'e\-pen\-dant de $l$. 
\endproclaim

\demo{Preuve}
Soient $l\not = l'$. Il suffit de composer les deux isomorphismes naturels de la preuve du th\'eor\`eme 15.2 pour obtenir un isomorphisme entre $K_l$ et $K_{l'}$.
$\square$
\enddemo

\proclaim{Corollaire 15.4}
Dans le cas particulier o\`u $T$ est une structure CR polaris\'ee de codimension z\'ero, l'espace $K_l$ (et l'espace $K_l^0$ qui lui est \'egal) s'identifient \`a l'espace de Kuranishi de la vari\'et\'e complexe compacte $(X,T)$.
\endproclaim

Ainsi les th\'eor\`emes 10.1 et 13.1 sont bien des g\'en\'eralisations du th\'eor\`eme de Kuranishi.

\vfill
\eject
\head
III. Vari\'et\'es CR $G$-polaris\'ees.
\endhead
\vskip1cm

Dans toute cette partie, $G$ est un groupe de Lie r\'eel connexe de dimension $d$. On notera $\frak G$ son alg\`ebre de Lie et $\exp_G$ son application exponentielle.

\vskip.5cm
\subhead
{\bf 16. Vari\'et\'es CR $\boldkey G$-polaris\'ees}
\endsubhead
\numcount=0

\definition{D\'efinition}
Une vari\'et\'e lisse compacte $X$ de dimension $2p+d$ est une {\it vari\'et\'e CR $G$-polaris\'ee} si elle est munie
\medskip
\noindent (i) d'une action localement libre de $G$.

\noindent (ii) d'une distribution CR $E$ de dimension $2p$ transverse aux orbites de $G$ et res\-pec\-t\'ee par l'action.
\enddefinition

Lorsque $d$ vaut $0$, le groupe $G$ est r\'eduit \`a un point et une vari\'et\'e CR $G$-polaris\'ee n'est rien d'autre qu'une vari\'et\'e complexe.

\medskip
Lorsque $p$ vaut $0$, le groupe $G$ a m\^eme dimension que $X$, si bien qu'une vari\'et\'e CR $G$-polaris\'ee est alors une vari\'et\'e homog\`ene (r\'eelle) $G/\Gamma$ avec $\Gamma$ r\'eseau de $G$.
\medskip
Ainsi cette notion est en quelque sorte une interpolation entre vari\'et\'e complexe et vari\'et\'e homog\`ene. Nous verrons cependant dans les sections suivantes qu'elle est en fait plus proche des vari\'et\'es complexes. 
\medskip
L'action localement libre de $G$ induit un feuilletage r\'eel $\Cal F$ sur $X$ et on a une d\'ecomposition du fibr\'e tangent
$$
TX=T\Cal F\oplus E.
\leqno (16.\numerote)
$$
De plus, la $G$-invariance du point (ii) signifie que 
$$
(g\cdot)_*E^{0,1}=E^{0,1}\quad\text{ pour tout }g\in G
\leqno (16.\numerote)
$$
o\`u $g\cdot$ d\'esigne l'action de $g$ sur $X$.

\proclaim{Proposition 16.1}
Une vari\'et\'e CR $G$-polaris\'ee est une vari\'et\'e CR polaris\'ee par $N=T\Cal F$.
\endproclaim

\demo{Preuve}
Le pseudo-groupe d'holonomie de $\Cal F$ est constitu\'e d'\'el\'ements de $G$ agissant sur le fibr\'e normal. Il r\'esulte donc de (16.\lastnum[0]) que $\Cal F$ est transversalement holomorphe et que $E$ est une r\'ealisation int\'egrable de sa structure holomorphe transverse. On conclut par le corollaire 3.2.
$\square$
\enddemo 

On peut construire sur $X$ des cartes feuillet\'ees de la mani\`ere classique suivante. Etant donn\'e $x\in X$ on se donne une section locale transverse au feuilletage 
$$
i\ :\ (\Bbb C^p,0)\longrightarrow (X,x)
\leqno (16.\numerote)
$$
et on pose
$$
(z,t)\in U\times\frak G\longmapsto (\exp_G t)\cdot i(z)\in X.
\leqno (16.\numerote)
$$
Nous appellerons de telles cartes {\it $G$-adapt\'ees}.
\medskip
Finissons cette section en indiquant quels exemples de la section 4 sont $G$-polaris\'es.

\example{Vari\'et\'es sasakiennes}
C'est l'exemple standard de structure CR $G$-polaris\'ee. Ici, le groupe $G$ est $\Bbb R$ qui agit sur $X$ via le flot du champ de Reeb d\'efini en (4.3). Nous avons d\'ej\`a vu en section 4 que ce flot pr\'eserve la structure CR.
\endexample

\example{Structures CR de codimension r\'eelle une invariantes sous l'action d'un flot transverse}
Elles sont, elles aussi et pour les m\^emes raisons, automatiquement $\Bbb R$-polaris\'ees, cf. section 4. Ces structures apparaissent dans la litt\'erature sasakienne sous le nom de {\it structures presque de contact normales}, cf. \cite{Bl}, \cite{B-G}.
\endexample

\example{Feuilletages transversalement holomorphes de codimension complexe u\-ne}
Nous avons vu en section 4 qu'un tel feuilletage donnait automatiquement une structure CR polaris\'ee sur $X$, en choisissant une distribution $E$ transverse. Mais cette fois une telle structure n'est pas forc\'ement $G$-polaris\'ee. D'une part, il faut que ce feuilletage provienne d'une action localement libre d'un groupe $G$. D'autre part, m\^eme en supposant une telle action, il faut encore que $G$ pr\'eserve la structure CR sur $E$, condition qui va d\'ependre fortement du choix de cette distribution transverse. 
\medskip
Dans le cas particulier des flots transversalement holomorphes sur les 3-vari\'et\'es, on voit dans la classification de Brunella-Ghys que tous les exemples poss\`edent, par construction,  une structure CR transverse naturelle. On peut pr\'eciser lesquelles sont $\Bbb R$-polaris\'ees. Ce travail est fait dans \cite{Ma} (voir aussi \cite{Ge} qui n'utilise pas la classification de Brunella-Ghys, mais directement le fait qu'une structure $\Bbb R$-polaris\'ee admet une complexification, cf. la section suivante). On trouve que tous les exemples sont $\Bbb R$-polaris\'es sauf les feuilletages stables des suspensions d'un diff\'eomorphisme hyperbolique de $T^2$ et les feuilletages transversalement affines de $\Bbb S^2\times\Bbb S^1$.  
\endexample

\example{Suspensions \`a base complexe}
Elles deviennent $G$-polaris\'ees lorsque la fibre $F$ est un groupe de Lie $G$, forc\'ement compact, et que la repr\'esentation $\rho$ est \`a valeurs dans les translations du groupe.
\endexample

\vskip.5cm
\subhead
{\bf 17. Le $\boldkey G$-c\^one et les $\boldkey G$-fibr\'es plats}
\endsubhead
\numcount=0

Ainsi que nous l'avons rappel\'e en section 4, une vari\'et\'e riemannienne est sasakienne si son c\^one riemannien (4.1) admet une structure k\"ahl\'erienne invariante par dilatations. Dans le cas ab\'elien, les vari\'et\'es CR $G$-polaris\'ees admettent une ca\-rac\-t\'e\-ri\-sation analogue. Pour l'\'enoncer, nous avons besoin de quelques d\'efinitions.

\definition{D\'efinition}
Soit $X$ une vari\'et\'e lisse. On appelle {\it $G$-c\^one} de $X$ le fibr\'e trivial $G\times X\to X$. On le notera $\Cal C_G(X)$.
\enddefinition

Pour faire court, on appellera $G$-fibr\'e un fibr\'e principal de fibre et de groupe structural $G$.
On rappelle qu'un $G$-fibr\'e est plat s'il poss\`ede une connexion de courbure nulle. Un $G$-fibr\'e plat sur une base simplement connexe est trivial \cite{No, Ch. II, \S 8}. Tout $G$-fibr\'e plat sur $X$ est donc quotient du $G$-c\^one de son rev\^etement universel $\tilde X$.

\definition{D\'efinition}
Soit $\pi : P\to X$ un $G$-fibr\'e plat sur une vari\'et\'e lisse $X$ et soit $H\subset TP$ une connexion de courbure nulle. Soit $J$ une structure complexe sur $P$. On dira que
\medskip
\noindent (i) $J$ est {\it invariante par translations} si $J$ est invariante par l'action naturelle de $G$ sur les fibres de $P$.

\noindent (ii) $J$ est {\it orthogonale aux fibres} si   
$$
v\in\text{Ker }d\pi\Longrightarrow Jv\in H.
\leqno (17.\numerote)
$$
\enddefinition

Observons que, si $P$ poss\`ede une structure complexe orthogonale aux fibres, les fibres de $P$ sont {\it totalement r\'eelles} dans $P$, i.e. pour tout $x\in X$, le fibr\'e tangent complexe de la fibre en $x$ est r\'eduit \`a z\'ero:
$$
T\pi^{-1}(x)\cap JT\pi^{-1}(x)=\{0\}.
\leqno (17.\numerote)
$$

Lorsque $P$ est le $G$-c\^one de $X$, la condition d'orthogonalit\'e aux fibres se r\'e\'ecrit plus simplement. Ecrivant abusivement
$$
T\Cal C_G(X)=TX\oplus \frak G
\leqno (17.\numerote)
$$
on demande que 
$$
J(\frak G)\subset TX.
\leqno (17.\numerote)
$$

La propri\'et\'e fondamentale des vari\'et\'es $G$-polaris\'ees ab\'eliennes est la suivante.

\proclaim{Th\'eor\`eme 17.1}
Soit $X$ vari\'et\'e lisse compacte. On identifie $X$ \`a l'hypersurface $X\times \{e\}$ de son $G$-c\^one. Soit $G$ un groupe de Lie connexe.
\medskip
\noindent i) Supposons $G$ ab\'elien et supposons que $X$ admette une structure CR $G$-polaris\'ee $(E,J)$. Alors $J$ s'\'etend en une structure complexe invariante par translations et orthogonale aux fibres sur $\Cal C_G(X)$.
\medskip
\noindent ii) R\'eciproquement, supposons que $\Cal C_G(X)$ admette une structure complexe invariante par translations et orthogonale aux fibres. Alors la structure CR induite sur $X$ est $G$-polaris\'ee et $G$ est ab\'elien.
\endproclaim

Plus g\'en\'eralement, on a le

\proclaim{Th\'eor\`eme 17.2}
Soit $X$ une vari\'et\'e CR $G$-polaris\'ee. On suppose $G$ ab\'elien. Soit $P$ un $G$-fibr\'e plat sur $X$ et $H$ une connexion plate. 

Alors il existe une structure complexe $J$ sur $P$ munie de $H$ invariante par translations et orthogonale aux fibres. De surcro\^{\i}t, la projection naturelle
$$
\pi_*\ :\ H\longrightarrow X
\leqno (17.\numerote)
$$
projette $J$ sur la structure CR de $X$.
\endproclaim

et le 

\proclaim{Th\'eor\`eme 17.3}
Soit $X$ vari\'et\'e lisse compacte. Supposons qu'il existe un $G$-fibr\'e plat $P$ sur $X$ de connexion plate $H$ admettant une structure complexe $J$ invariante par translations et orthogonale aux fibres. 

Alors la projection $\pi_*J$ (cf. {\rm (17.\lastnum[0])}) d\'efinit une structure $G$-polaris\'ee sur $X$. De plus, $G$ est ab\'elien.
\endproclaim

\demo{Preuve des th\'eor\`emes 17.1, 17.2 et 17.3}
Le th\'eor\`eme 17.1 est une application directe des th\'eor\`emes 17.2 (pour la partie i)) et 17.3 (pour la partie ii)) en utilisant la connexion triviale (17.\lastnum[-2]).
\medskip
Commen\c cons donc par prouver le th\'eor\`eme 17.2. Soit $P$ un $G$-fibr\'e plat et $H$ une connexion plate sur $G$. Nous allons construire deux types de champs fondamentaux. Soit $\xi\in\frak G$. On pose d'une part
$$
(x,g)\in P\longmapsto \xi_*(x,g):=H_{(x,g)}\left (\dfrac{d}{ds}\vert_{s=0}(\exp_G(s\xi)\cdot x)\right )\in T_{(x,g)}P
\leqno (17.\numerote)
$$
o\`u  $H_{(x,g)}(w)$ est le relev\'e horizontal en $(x,g)$ du vecteur $w\in T_xX$; et d'autre part,
$$
(x,g)\in P\longmapsto \xi^*(x,g):=\dfrac{d}{ds}\vert_{s=0}\bigg (x, \exp_G(s\xi)\cdot g\bigg )\in T_{(x,g)}P.
\leqno (17.\numerote)
$$

Par d\'efinition, $\xi_*$ est un champ horizontal, tandis que $\xi^*$ est un champ vertical. 
\medskip
Appelons comme d'habitude $J$ la structure CR. On \'etend $J$ \`a $P$ de la mani\`ere suivante.   On pose d'une part
$$
J_{(x,g)}v:=H_{(x,g)}\left ( J_x(\pi_*v)\right )
\leqno (17.\numerote)
$$
pour $v\in H_{(x,g)}$ v\'erifiant $\pi_*v\in E$ ; et d'autre part
$$
J\xi_*=\xi^*\qquad\text{ et }\qquad J\xi^*=-\xi_*
\leqno (17.\numerote)
$$
pour $\xi\in \frak G$. On obtient ainsi un op\'erateur presque complexe sur $P$. 
\medskip
Montrons qu'il est int\'egrable. Pour cela, notons que
$$
TP^{0,1}=H^{0,1}\oplus \text{Vect}_{\Bbb C}\{\xi_*+i\xi^*\quad\vert\quad\xi\in\frak G\}
\leqno (17.\numerote)
$$
o\`u $H^{0,1}$ est le relev\'e horizontal de $E^{0,1}$.

\medskip 
Observons que $H^{0,1}$ est involutif puisque $E^{0,1}$ l'est et que la connexion est plate. Pla\c cons-nous localement et utilisons les champs $e_i$. Soit $\xi\in\frak G$. Comme l'action de $G$ pr\'eserve $E^{0,1}$, on a la propri\'et\'e infinit\'esimale (cf. (4.4) et (4.6)) 
$$
[e_j,\pi_*\xi_*]\in E^{0,1}
\leqno (17.\numerote)
$$
et comme par d\'efinition $\pi_*\xi_*$ est dans $N$, on a de plus
$$
[e_j,\pi_*\xi_*]\in E^{0,1}\cap N=\{0\}.
\leqno (17.\numerote)
$$
d'o\`u 
$$
[He_j,\xi_*]=0.
\leqno (17.\numerote)
$$
Maintenant, la d\'efinition (17.\lastnum[-6]) montre que $\xi^*$ ne d\'epend pas des coordonn\'ees en $X$. En fait, c'est imm\'ediat si le fibr\'e est trivial ; et c'est une propri\'et\'e qui passe aux quotients par un groupe discret, donc c'est vrai pour tout fibr\'e plat. Ainsi
$$
[He_j,\xi^*]=0.
\leqno (17.\numerote)
$$
Enfin, on a, pour tout $\xi\in \frak G$, et tout $\eta\in\frak G$,
$$
[\xi_*,\eta^*]=0.
\leqno (17.\numerote)
$$
On conclut de (17.\lastnum[0]), (17.\lastnum[-1]) et (17.\lastnum[-2]) et de la commutativit\'e de $G$ que $TP^{0,1}$ est involutif.
\medskip
Les d\'efinitions (17.\lastnum[-9]) et (17.\lastnum[-8]) impliquent imm\'ediatement que la structure complexe ainsi d\'efinie est orthogonale aux fibres. Enfin, on d\'eduit de (17.\lastnum[-1]) et (17.\lastnum[0]) et de la commutativit\'e de $G$ que
$$
[\xi^*, TP^{0,1}]\subset TP^{0,1}
\leqno (17.\numerote)
$$
pour tout $\xi\in\frak G$, donc que la structure $J$ est invariante par translations. Le th\'eor\`eme 17.2 est d\'emontr\'e.
\medskip
Montrons maintenant le th\'eor\`eme 17.3. Soit $X$ vari\'et\'e lisse compacte. Supposons qu'il existe un $G$-fibr\'e plat $P$ sur $X$ de connexion plate $H$ admettant une structure complexe $J$ invariante par translations et orthogonale aux fibres. 
\medskip
Pour tout $\xi\in\frak G$, on d\'efinit $\xi^*$ comme en (17.\lastnum[-9]), et on pose
$$
\xi_*:=J\xi^*.
\leqno (17.\numerote)
$$
La propri\'et\'e d'orthogonalit\'e aux fibres signifie que $\xi_*$ est un champ horizontal. L'invariance par translations implique que les $\xi_*$ sont $G$-invariants. On en d\'eduit que
$$
[\xi_*,\eta^*]=0
\leqno (17.\numerote)
$$
puis, en \'ecrivant la nullit\'e du tenseur de Niejenhuis, 
$$
[\xi^*,\eta^*]=-[\xi_*,\eta_*].
\leqno (17.\numerote)
$$
Mais, le premier crochet est vertical, et par platitude le deuxi\`eme horizontal, donc ils sont en fait nuls. Ceci prouve que $G$ est ab\'elien, mais aussi que l'ensemble des champs $\pi_*\xi_*$ lorsque $\xi$ varie dans $\frak G$ forme une sous-alg\`ebre de Lie des champs de vecteurs de $X$ isomorphe \`a $\frak G$. Autrement dit, on d\'efinit ainsi une action de $G$ sur $X$. Comme ces champs ne s'annulent pas, cette action est localement libre.
\medskip
En outre, on v\'erifie ais\'ement que, si l'on pose
$$
\Bbb CH:=H\cap JH
\leqno (17.\numerote)
$$
alors d'une part on a
$$
\Bbb CH=\{v+JV(Jv)\quad\vert\quad v\in H\}
\leqno (17.\numerote)
$$
pour $V$ projection verticale de $H\oplus\frak G$ sur $\frak G$ ; et d'autre part, on en d\'eduit que
$$
H=\Bbb CH\oplus J\frak G.
\leqno (17.\numerote)
$$
On peut donc d\'efinir 
$$
H^{0,1}:=\{v+iJv\quad\vert\quad v\in\Bbb CH\}
\leqno (17.\numerote)
$$ 
et la projection
$$
E^{0,1}=\pi_*H^{0,1}
\leqno (17.\numerote)
$$
munit $X$ d'une structure CR transverse au feuilletage d\'efini par l'action de $G$.
\medskip
La propri\'et\'e d'invariance par translations signifie que cette structure est invariante par l'action de $G$. Il s'agit donc bien d'une structure CR $G$-polaris\'ee et le th\'eor\`eme 17.3 est d\'emontr\'e.
$\square$
\enddemo

 \vskip.5cm
\subhead
{\bf 18. Vari\'et\'es CR $\boldkey G$-polaris\'ees et vari\'et\'es compactes complexe}
\endsubhead
\numcount=0

Le but de cette section est d'associer \`a une vari\'et\'e CR $G$-polaris\'ee des vari\'et\'es compactes complexes.
\medskip
Le premier r\'esultat est le

\proclaim{Th\'eor\`eme 18.1}
Soit $X$ une vari\'et\'e CR $G$-polaris\'ee. On suppose $G$ compact de dimension paire. Alors $X$ est une vari\'et\'e compacte complexe, la structure CR sur $E$ est induite de $X$ et le feuilletage $\Cal F$ est holomorphe.
\endproclaim

\demo{Preuve}
Par le fameux th\'eor\`eme de Samelson \cite{Sa}, $G$ poss\`ede une structure complexe $J$ invariante par les translations \`a gauche. Comme l'action de $G$ sur $X$ est localement libre, les champs fondamentaux de l'action
$$
\xi_*(x)=\left (\dfrac{d}{ds}_{\vert s=0} (\exp (s\xi)\cdot x) \right )
\quad\text{ pour } \xi\in\frak G,\ x\in X
\leqno (18.\numerote)
$$
ne s'annulent pas si bien que l'on peut d\'efinir une structure presque complexe le long des feuilles de $\Cal F$ par
$$
J_{\Cal F}\xi_*:=(J\xi)_*.
\leqno (18.\numerote)
$$
Comme par ailleurs,
$$
[\xi_*,\eta_*]=[\xi,\eta]_*
\leqno (18.\numerote)
$$
cette structure est int\'egrable et on munit ainsi $\Cal F$ d'une structure complexe. Maintenant, on a
$$
TX=E\oplus T\Cal F
\leqno (18.\numerote)
$$
si bien qu'on peut d\'efinir
$$
J_X:=J\oplus J_\Cal F
\leqno (18.\numerote)
$$
sur $X$. Il reste \`a v\'erifier que cette structure presque complexe est int\'egrable. Mais localement les champs $(0,1)$ de $E$ sont engendr\'es par les $e_i$ et les champs $(0,1)$ de $\Cal F$ par les 
$$
\xi_*+i(J\xi)_*\qquad\qquad \xi\in\frak G.
\leqno (18.\numerote)
$$
Or, les premiers sont stables par crochet par int\'egrabilit\'e de la structure CR de $E$, et on vient de montrer que les seconds le sont. Enfin, les crochets
$$
[e_i, \xi_*+i(J\xi)_*]
\leqno (18.\numerote)
$$
restent $(0,1)$ par $G$-invariance de $E^{0,1}$ (cf. (17.11)). Par d\'efinition cette structure est feuillet\'ee et $\Cal F$ est donc holomorphe. Enfin, par d\'efinition, elle induit la structure $J$ sur $E$.
$\square$
\enddemo

\proclaim{Corollaire 18.2}
Si $G$ est compact, alors $X$ ou $X\times \Bbb S^1$ admet une structure complexe telle que la structure CR sur $E$ soit induite de $X$ et telle que le $G$-feuilletage $\Cal F$ de $X$, respectivement le $G\times\Bbb S^1$-feuilletage de $X\times\Bbb S^1$ soit holomorphe.
\endproclaim

\demo{Preuve}
Pour $G$ de dimension paire, c'est exactement le th\'eor\`eme 18.1. Si $G$ est de dimension impaire, alors $G\times \Bbb S^1$ agit localement librement sur $X\times\Bbb S^1$. Il suffit maintenant de transporter la structure CR de $X$ sur $X\times\Bbb S^1$ en faisant agit $G\times \Bbb S^1$ sur $X\times\{e\}\subset X\times \Bbb S^1$ pour obtenir une $G\times\Bbb S^1$-polarisation de $X\times\Bbb S^1$. On conclut de nouveau par le th\'eor\`eme 18.1.
$\square$
\enddemo

Remarquons que la preuve du th\'eor\`eme 18.1 n'utilise pas la compacit\'e et s'\'etend au cas o\`u $G$ est non compact mais admet une structure de groupe de Lie complexe, ou tout au moins une structure complexe invariante par translations \`a gauche. On a donc

\proclaim{Corollaire 18.3}
Soit $X$ une vari\'et\'e CR $G$-polaris\'ee. On suppose que $G$ admet une structure complexe invariante par translations \`a gauche. Notons $G_\Bbb C$ cette vari\'et\'e complexe. Alors, 
\medskip
\noindent (i) La vari\'et\'e $X$ est une vari\'et\'e compacte complexe $X_\Bbb C$, la structure CR sur $E$ est induite de $X_\Bbb C$ et le feuilletage $\Cal F$ est holomorphe.

\noindent (ii) Soit $P$ un $G$-fibr\'e principal plat sur $X$. Alors il existe une structure complexe sur $P$ telle que la projection naturelle sur $X$ soit un fibr\'e holomorphe localement trivial de base $X_\Bbb C$, de fibre $G_\Bbb C$ et de groupe structural $G$.
\endproclaim

\remark{Remarque}
Pour que le (ii) soit valable, il faut supposer que le groupe structural de $P$ agisse \`a gauche sur les fibres de $P$.
\endremark

\demo{Preuve}
Il nous reste uniquement \`a montrer le (ii). Recouvrons $P$ de cartes locales de fibr\'e, i.e. localement diff\'eomorphes \`a $X\times G$. Comme le fibr\'e $P$ est plat, on peut supposer les changements de cartes du type
$$
(x,g)\longmapsto (\phi(x), \psi(x)\cdot g)
\leqno (18.\numerote)
$$
avec $\psi$ localement constante. On munit $P$ localement de la structure complexe produit des cartes. Les changements de cartes (18.\lastnum[0]) pr\'eservent cette structure complexe, puisqu'ils sont constitu\'es d'une part de changements de cartes holomorphes de $X$, d'autre part de translations \`a gauche de $G$ localement constantes.
$\square$
\enddemo

 En particulier, le corollaire 18.3 s'applique au cas $G$ ab\'elien de dimension paire et permet de munir tout $G$-fibr\'e plat sur $X$ d'une structure complexe. Il est particuli\`erement int\'eressant de comparer cette structure \`a celle donn\'ee par le th\'eor\`eme 17.2. 

\proclaim{Proposition 18.4}
Soit $X$ une vari\'et\'e CR $\Bbb R^m$-polaris\'ee. 
Supposons $m$ pair. Soit $P$ un $G$-fibr\'e plat sur $X$. Alors, la structure complexe de $P$ donn\'ee par le th\'eor\`eme 17.2 n'est pas biholomorphe \`a celle donn\'ee par le corollaire 18.3.
\endproclaim

\demo{Preuve}
On constate que la structure complexe du th\'eor\`eme 17.1, i) est orthogonale aux fibres de $\Cal C_{\Bbb R^m}(X)\to X$; alors que celle du corollaire 18.3 est une structure de fibr\'e holomorphe principal sur cette projection.
$\square$
\enddemo  

\remark{Remarque}
Lorsque $G$ vaut $\Bbb R$, il existe \'egalement deux structures complexes sur le $\Bbb R$-c\^one de $X$ : celle donn\'ee par le th\'eor\`eme 17.1, i) et celle donn\'ee par le corollaire 18.2 (en passant au rev\^etement $X\times\Bbb R\to X\times\Bbb S^1$). Mais une lecture attentive des preuves de ces deux r\'esultats montre qu'elles sont biholomorphes.
\endremark
\medskip

Dans certains cas, on peut associer \`a $X$ une autre vari\' et\' e compacte complexe. La d\'efinition suivante est consistante avec \cite{Sp} et \cite{B-G}.

\definition{D\'efinition}
Une vari\'et\'e CR $G$-polaris\'ee $X$ est dite {\it r\'eguli\`ere} lorsque l'action localement libre de $G$ sur $X$ provient d'une action libre d'un groupe compact $H=G/\Gamma$ pour $\Gamma$ r\'eseau de $G$. On appellera {\it quotient r\'egulier de} $G$ un tel groupe. 

Elle est dite {\it quasi-r\'eguli\`ere} si l'action provient d'une action avec stabilisateurs finis d'un groupe compact $H=G/\Gamma$ pour $\Gamma$ r\'eseau de $G$. On appellera {\it quotient quasi-r\'egulier de} $G$ un tel groupe. 

Enfin, elle est dite {\it irr\'eguli\`ere} si elle n'est pas quasi-r\'eguli\`ere.
\enddefinition

\proclaim{Proposition 18.5}
Soit $X$ une vari\' et\' e CR $G$-polaris\' ee r\'eguli\`ere. Soit $H$ un quotient r\'egulier de $G$. Alors,
\medskip
\noindent (i) le quotient $X/H$ est une vari\'et\'e compacte complexe $Y$.

\noindent (ii) Si de plus $G$ est de dimension paire, alors la projection $X\to Y$ est un fibr\'e holomorphe $H$-principal.
\endproclaim

\demo{Preuve}
L'action de $H$ est libre et propre par d\'efinition, donc le quotient $Y$ est une vari\' et\'e lisse. Par ailleurs, par d\'efinition, $H$ engendre la m\^eme action que $G$, et donc le feuilletage qu'il induit sur $X$ est exactement le feuilletage $\Cal F$. Le quotient $Y$ s'identifie donc \`a l'espace des feuilles de $\Cal F$, et la structure CR transverse $G$-invariante descend en une structure complexe.
\medskip
Pour le (ii), on applique le th\'eor\`eme 18.1 \`a $H$ et on note que la structure complexe ainsi obtenue sur $X$ est pr\'eserv\'ee par l'action de $H$. On conclut par le th\'eor\`eme de Holmann.
$\square$
\enddemo
 
Ceci motive la

\definition{D\'efinition}
On appellera {\it vari\'et\'e quotient de $X$}  une vari\'et\'e compacte complexe construite de cette mani\`ere.
\enddefinition

\remark{Remarque}
Lorsque $X$ est quasi-r\' eguli\`ere, on obtient un quotient $Y$ qui est une orbifold. On peut facilement adapter tous les r\'esultats qui suivent au cas quasi-r\' egulier.
\endremark

\vskip.5cm
\subhead
{\bf 19. Le cas ab\'elien}
\endsubhead
\numcount=0

Supposons $G$ ab\'elien.
Soit $X$ une vari\'et\'e CR $G$-polaris\'ee. Le th\'eor\`eme 17.2 entra\^{\i}ne que tout $G$-fibr\'e plat $P$ sur $X$ admet une structure complexe invariante par translations. Soit $\Gamma$ un r\'eseau cocompact. Le quotient
$$
Z=P/\Gamma
\leqno (19.\numerote)
$$
o\`u $\Gamma$ agit sur les fibres de $P$ via son inclusion dans $G$, est une vari\'et\'e compacte complexe.

\definition{D\'efinition}
On appellera {\it vari\'et\'e associ\'ee \`a $P$}  une vari\'et\'e compacte complexe construite de cette mani\`ere.
\enddefinition

Par d\'efinition, une vari\'et\'e associ\'ee \`a $P$ est un fibr\'e {\it r\'eel} $G/\Gamma$-principal sur $X$. {\it Mais ce n'est pas un fibr\'e holomorphe} sur $X$, puisque les fibres $G/\Gamma$ sont totalement r\'eelles dans $P$ (cf. Proposition 18.4).
\medskip
Pla\c cons-nous dans le cas o\`u l'action de $G$ est r\'eguli\`ere. Soit $H$ le groupe quotient associ\'e. Nous pouvons donc d\'efinir la vari\'et\'e quotient $Y$.
Nous souhaitons maintenant \' eclaircir la relation entre $P$ et $Y$. Le feuilletage de $X$ peut \^etre remont\'e \`a $P$ en utilisant la connexion plate. Ce feuilletage de $P$ est lui aussi donn\'e par une action dont les champs fondamentaux sont les images par $J$ des champs fondamentaux de l'action de $G$ sur les fibres de P (cf. (17.6), (17.7)). Appelons {\it horizontale} cette action et {\it verticale} l'action sur les fibres de $P$.
\medskip
On peut donc faire agir $G\times G$ sur $P$, le premier facteur agissant horizontalement et le deuxi\`eme verticalement. Appelons {\it compl\`ete} cette action de $G\times G$ sur $P$. 
\medskip
Notons que le groupe ab\'elien $G\times G$ a une structure de groupe de Lie complexe obtenue en d\'ecr\'etant que dans la d\'ecomposition de son alg\`ebre de Lie
$$
\frak K=\frak G\oplus \frak G
\leqno (19.\numerote)
$$
le premier terme est la partie r\'eelle et le second la partie imaginaire. On notera $(G\times G)_{\Bbb C}$ ce groupe de Lie complexe.

Observons que la structure complexe de $(G\times G)_\Bbb C$ descend en une structure de groupe de Lie complexe sur $G\times H$ et sur $G/\Gamma \times H$. On les notera $(G\times H)_\Bbb C$ et $(G/\Gamma \times H)_\Bbb C$.

\proclaim{Proposition 19.2}
Soit $X$ une vari\' et\' e CR $G$-polaris\' ee r\'eguli\`ere avec $G$ ab\'elien. Soit $H$ un quotient r\' egulier de $G$. Soit enfin $P$ un $G$-fibr\' e plat sur $X$. On munit $P$ de la structure complexe donn\'ee par le th\'eor\`eme 17.2. Alors
\medskip
\noindent (i) L'action compl\`ete de $(G\times G)_{\Bbb C}$ sur la vari\'et\'e complexe $P$ est holomorphe.

\noindent (ii) La vari\' et\' e quotient $Y$ est le quotient de $P$ par l'action compl\`ete de $G\times G$.

\noindent (iii) La vari\'et\'e $P$ est un fibr\' e holomorphe principal sur $Y$ de groupe $(G\times H)_\Bbb C$.

\noindent (iv) On a un diagramme commutatif
$$
\CD
P \aro>\Gamma >> Z\cr
\aro V(G\times H)_{\Bbb C} VV  \aro VV (G/\Gamma\times H)_{\Bbb C} V\cr
Y\aro >>Id> Y
\endCD
\leqno (19.\numerote)
$$
dont les fl\`eches verticales sont des fibr\'es holomorphes principaux.
\endproclaim

\demo{Preuve}
On observe que dans la preuve du th\'eor\`eme 17.2, la structure complexe des orbites de l'action compl\`ete co\"{\i}ncide avec celle de $(G\times G)_{\Bbb C}$. De plus,
chaque \'el\'ement de $G$ agissant verticalement agit holomorphiquement, par invariance par translations de la structure complexe de $P$. De m\^eme, comme l'action horizontale de $G$ sur $P$ est infinit\'esimalement $J$ fois la verticale, chaque \'el\' ement agissant horizontalement agit holomorphiquement. Ceci prouve le (i).
\medskip
Prenons le quotient de $P$ par l'action compl\`ete en deux temps. D'une part, par d\'efinition, le quotient de $P$ par l'action verticale est $X$ et l'action horizontale descend en l'action de $G$ sur $X$. D'autre part, par r\'egularit\'e, cette action se r\'eduit \`a une action effective de $H$ dont le quotient est $Y$. On obtient ainsi un diagramme commutatif (18.\lastnum[0]), mais pour l'instant sans savoir qu'il s'agit d'un diagramme de fibr\' es.
\medskip
L'action compl\`ete de $G\times H$ sur $P$ est une action libre, propre et d'apr\`es le (i), holomorphe avec orbites qui s'identifie \`a $(G\times H)_\Bbb C$. Il est alors facile de terminer la preuve.
$\square$
\enddemo
  
Finissons cette section avec l'exemple suivant.

\example{Vari\'et\'es sasakiennes}
Faisons un petit tour d'horizon des notions pr\' ec\'edentes lorsque $X$ est une vari\' et\'e sasakienne. On a alors $G=\Bbb R$ et le $\Bbb R$-c\^one s'identifie au c\^one riemannien (diff\'erentiablement, le $\Bbb R$-c\^one n'\' etant pas muni d'une m\'etrique particuli\`ere) par
$$
(x,t)\in \Cal C_{\Bbb R} (X)\longmapsto (x,\exp t)\in \Cal C(X).
\leqno (19.\numerote)
$$
Il n'y a pas d'autre fibr\'es plats, tout $\Bbb R$-fibr\' e plat \' etant automatiquement trivial.
\medskip

Le cas r\'egulier (cf. \cite{Sp}) correspond au cas o\`u le flot provient d'une action libre de $\Bbb S^1$. Il est alors bien connu que la vari\'et\'e quotient $Y$ est une vari\'et\'e projective et que $X$ est un fibr\' e unitaire d'un fibr\'e en droites sur $Y$. Enfin le c\^one riemannien muni de sa structure complexe est le $\Bbb C^*$-fibr\'e principal associ\'e sur $Y$ (cf. \cite{Sp} pour des \'enonc\'es plus pr\'ecis).
\endexample

 \vskip.5cm
\subhead
{\bf 20. Vari\'et\'es CR $\Bbb R^m$-polaris\'ees et vari\'et\'es LVM}
\endsubhead
\numcount=0

Le but de cette section est de relier les vari\'et\'es CR $\Bbb R^m$-polaris\'ees et les vari\'et\'es LVM. On renvoie \`a \cite{LdM-V}, \cite{Me1} et \cite{M-V} pour une pr\'esentation de ces vari\'et\'es compactes complexes non k\"ahl\'eriennes. Nous utiliserons directement les notations et r\'esultats de \cite{M-V}.
\medskip
Pour prendre la mesure du th\'eor\`eme 20.2 \'enonc\'e ci-dessous, remarquons qu'on d\'eduit facilement de la section 17 le r\'esultat suivant.

\proclaim{Proposition 20.1}
Soit $Y$ une vari\'et\'e compacte complexe. Soit $m$ un entier naturel et soient $L_1,\hdots, L_m$ des fibr\'es en droites sur $Y$. Munissons-les de m\'etriques riemanniennes et appelons $C_1,\hdots ,C_m$ les fibr\'es en cercles unitaires associ\'es. Appelons $X$ la vari\'et\'e $C_1\oplus\hdots\oplus C_m$.

Alors $X$ est une vari\'et\'e CR $\Bbb R^m$-polaris\'ee r\' eguli\`ere pour $H=(\Bbb S^1)^m$. Son $\Bbb R^m$-c\^one est le $(\Bbb C^*)^m$ fibr\'e principal 
$$
P:=L_1\setminus\{0\}\oplus\hdots\oplus L_m\setminus\{0\}
\leqno (20.\numerote)
$$
et sa vari\'et\'e quotient est $Y$.
\endproclaim

\demo{Preuve}
Choisissons un flot engendrant chaque action circulaire sur $X$. On obtient ainsi une $\Bbb R^m$-action sur $X$ dont le quotient est $Y$.
On voit imm\'ediatement que $P$ est diff\'erentiablement le $\Bbb R^m$-c\^one de $X$ et que la structure complexe de $P$ est invariante par translations et orthogonale aux fibres. On conclut par le th\'eor\`eme 17.1, ii).
$\square$
\enddemo

On a m\^eme mieux. Un fibr\'e en cercles sur une vari\'et\'e complexe $Y$ est $\Bbb R$-polaris\'e si et seulement si c'est le fibr\'e unitaire d'un fibr\'e en droites sur $Y$. Il s'agit d'une cons\'equence directe de \cite{H-S}.
\medskip

Ceci montre qu'il est tr\`es facile de construire des vari\'et\'es CR $\Bbb R^n$-polaris\'ees r\'eguli\`eres. Comme dans le cas sasakien, il est beaucoup plus difficile de construire des exemples {\it irr\'eguliers}. Pour y arriver, nous allons utiliser les vari\'et\'es LVM. Rappelons que ces vari\'et\'es non k\"ahl\'eriennes poss\`edent toujours un feuilletage transversalement k\"ahl\'erien donn\'e par une action de $\Bbb C^m$. Qui plus est, sous une condition de rationalit\'e dite condition (K), ce feuilletage est \`a feuilles compactes (voir \cite{M-V, Theorem A}).
  
\definition{D\'efinition}
Soit $N$ une vari\'et\'e LVM. On dira qu'elle {\it v\'erifie la condition} (K0) si
\medskip
\noindent (i) L'action transversalement k\"ahl\'erienne de \cite{M-V, Theorem A} se r\'eduit \`a une action de $(\Bbb C^*)^m$.

\noindent (ii) L'action de $(\Bbb S^1)^m$ induite du point (i) est libre.
\enddefinition
Nous appellerons {\it action }(K0) l'action de $(\Bbb S^1)^m$ sur une vari\'et\'e LVM v\'erifiant la condition (K0).
\medskip
Remarquons que cette condition est loin d'\^etre vide. En effet, si $N$ v\'erifie la condition (K) de \cite{M-V}, le point (i) de la condition (K0) est v\'erifi\'e ; et si $N$ v\'erifie la condition de non-singularit\'e de \cite{M-V, Corollary B, (ii)}, elle satisfera de plus au point (ii). Le corollaire H de \cite{M-V} montre que de tels exemples existent au-dessus de n'importe quelle vari\'et\'e torique projective lisse, et ce avec un nombre de points indispensables arbitrairement grand (voir \cite{M-V, Definition 1.8}. Mais remarquons \'egalement, et c'est fondamental pour la suite, qu'il y a aussi de nombreux exemples o\`u $N$ v\'erifie la condition (K0) mais pas la condition (K), puisque cette derni\`ere implique que l'action transversalement k\"ahl\'erienne se r\'eduit \`a une action de $(\Bbb S^1)^{2m}$.
\medskip

On a le

\proclaim{Th\'eor\`eme 20.2}
Soit $N$ une vari\'et\'e LVM v\'erifiant la condition {\rm (K0)}. On suppose que $N$ poss\`ede au moins un point indispensable. Alors
\medskip
\noindent (i) Le quotient de $N$ par l'action {\rm (K0)} est une vari\'et\'e CR $\Bbb R^m$-polaris\'ee $X$.

\noindent (ii) La vari\'et\'e $X$ est quasi-r\'eguli\`ere si et seulement si $N$ v\'erifie la condition {\rm (K)}.
\endproclaim

Ce th\'eor\`eme donne ainsi de nombreux exemples de vari\'et\'es CR $\Bbb R^m$-polaris\'ees irr\'eguli\`eres.

\demo{Preuve}
Le quotient de $N$ par l'action libre et propre (K0) est une vari\'et\'e lisse $X$ munie d'une action de $\Bbb R^m$. En effet, l'action holomorphe commutative initiale de $\Bbb C^m$ sur $N$ est, de par la condition (K0), une action commutative de $(\Bbb C^*)^m$. Cette action se d\'ecompose naturellement en l'action (K0) d'une part et une action de $\Bbb R^m$ d'autre part. Ces deux actions commutant, la $\Bbb R^m$-action descend \`a $X$.

\proclaim{Lemme 20.3}
La $\Bbb R^m$-action sur $X$ est localement libre.
\endproclaim

\demo{Preuve du lemme 20.3}
Elle est induite de la $\Bbb C^m$-action transversalement k\"ahl\'erien\-ne qui est localement libre \cite{M-V,\S 2}.
$\square$
\enddemo

Le quotient $X$ est donc muni d'un feuilletage r\'eel de dimension $m$. La structure transversalement k\"ahl\'erienne de $N$ descend en une structure transversalement k\"ahl\'erien\-ne sur ce feuilletage. 
\medskip
Consid\'erons alors la $2$-forme transversalement k\"ahl\'erienne $\omega$ de $N$ (appel\'ee forme d'Euler canonique dans \cite{M-V}). Lorsqu'on identifie $N$ diff\'erentiablement avec la sous-vari\'et\'e r\'eelle $\Cal N$ de l'espace projectif complexe $\Bbb P^{n-1}$ (cf. \cite{M-V, formule (8)}), il s'agit simplement de la restriction \`a $\Cal N$ de la forme de Fubini-Study du projectif.
\medskip
L'hypoth\`ese de l'existence d'un point indispensable entra\^{\i}ne que $\Cal N$ vit en fait dans un ouvert affine du projectif, et donc que $\omega$ est exacte. Soit $\alpha$ une $1$-forme de $N$ dont la diff\'erentielle est $\omega$. On a maintenant le

\proclaim{Lemme 20.4}
On peut prendre $\alpha$ invariante par l'action {\rm (K0)}.
\endproclaim

\demo{Preuve du lemme 20.4}
En fait, si l'on suppose que $\Cal N$ est incluse dans la carte affine $(w_1,\hdots, w_{n-1})$ de $\Bbb P^{n-1}$ correspondant \`a $z_1=1$, elle s'\'ecrit
$$
\Cal N=\{w\in\Bbb C^{n-1}\quad\vert\quad \Lambda_1+\sum_{i=1}^{n-1}\Lambda_{i+1}\vert w_i\vert ^2=0\}.
\leqno (20.\numerote)
$$
Il r\'esulte de ce qui a \'et\'e rappel\'e de la forme $\omega$ qu'on peut prendre
$$
\alpha=\dfrac{\sum w_id\bar w_i-\bar w_i dw_i}{1+\sum\vert w_i\vert ^2}.
\leqno (20.\numerote)
$$
Par ailleurs, nous affirmons que l'action (K0) se fait par multiplication des coordonn\'ees par des nombres complexes de module $1$, donc pr\'eserve $\alpha$.
\medskip
Pour montrer cela, il faut revenir \`a la construction des vari\'et\'es LVM et de l'action transversalement k\"ahl\'erienne. Par d\'efinition, $N$ est la projectivisation de l'espace des feuilles de Siegel (i.e. ferm\'ees et ne contenant pas l'origine) du feuilletage de $\Bbb C^n$ donn\'e par l'action
$$
(T,z)\in \Bbb C^m\times\Bbb C^n\longmapsto \Bigl ( z_i\exp \langle \Lambda_i, T\rangle \Bigr )_{i=1}^n\in\Bbb C^n.
\leqno (20.\numerote)
$$
Nous renvoyons \`a \cite{M-V, \S 1} pour plus de d\'etails. En tenant compte du fait que $z_1$ est une coordonn\'ee indispensable, on peut r\'e\'ecrire $N$ comme le quotient d'un ouvert $U$ de $\Bbb C^{n-1}$ (correspondant aux feuilles de Siegel) par l'action
$$
(T,w)\in\Bbb C^m\times U\longmapsto A_T(w):=\Bigl ( w_i\exp \langle \Lambda_{i+1}-\Lambda_1, T\rangle \Bigr )_{i=1}^{n-1}\in U.
\leqno (20.\numerote)
$$
De plus, $N$ s'identifie au mod\`ele $C^\infty$ (20.\lastnum[-3]).
\medskip
Ensuite, l'action transversalement k\"ahl\'erienne est l'action induite sur $N$ par l'action 
$$
(S,w)\in\Bbb C^m\times U\longmapsto B_S(w):=\Bigl ( w_i\exp \langle \Re (\Lambda_{i+1}-\Lambda_1), T\rangle \Bigr )_{i=1}^{n-1}\in U
\leqno (20.\numerote)
$$
et l'action (K0) est la projection \`a $N$ de (20.\lastnum[0]) pour $S\in (\Bbb S^1)^m\subset \Bbb C^m$.
\medskip
Cette action ne pr\'eserve pas directement le mod\`ele (20.\lastnum[-4]). Toutefois, \'etant donn\'e $w\in\Cal N$ et $S\in (\Bbb S^1)^m$, un calcul direct montre que la feuille de l'action $A$ passant par $B_S(w)$ coupe $\Cal N$ en un point unique \'egal \`a 
$$
A_{-S}\circ B_S(w).
\leqno (20.\numerote)
$$
Autrement dit, l'action (K0) transpos\'ee au mod\`ele (20.\lastnum[-5]) s'\'ecrit
$$
(w,S)\in\Cal N\times (\Bbb S^1)^m\longmapsto \Bigl ( w_i\exp \langle \Im (\Lambda_{i+1}-\Lambda_1), iS\rangle \Bigr )_{i=1}^{n-1}\in\Cal N.
\leqno (20.\numerote)
$$
Ces formules prouvent ce que nous affirmions : l'action (K0) pr\'eserve la forme (20.\lastnum[-5]).
$\square$
\enddemo

La $1$-forme $\alpha$ descend donc en une $1$-forme $\beta$ sur $X$. Cette forme n'est, par contre, pas invariante par le feuilletage de $X$. Appelons $p$ la projection de $N$ sur $X$. La diff\'erentielle $p_*$ d\'efinit un morphisme surjectif
$$
p_*\ :\ \text{Ker }\alpha\subset TN\longrightarrow \text{Ker }\beta\subset TX.
\leqno (20.\numerote)
$$
Le point-clef est que ce morphisme envoie la distribution CR de codimension un de $N$
$$
(\bar E,\bar J)=\text{Ker }\alpha\cap i\text{Ker }\alpha
\leqno (20.\numerote)
$$
en la distribution CR de codimension $m$ de $X$
$$
(E,J)=p_*(\text{Ker }\alpha\cap i\text{Ker }\alpha)
\leqno (20.\numerote)
$$
invariante par l'action de $\Bbb R^m$, et transverse aux orbites de cette action. Ainsi $X$ est munie d'une structure CR $\Bbb R^m$-polaris\'ee.
Ceci montre le (i).
\medskip
Maintenant $X$ est quasi-r\'eguli\`ere si et seulement si la $\Bbb R^m$-action sur $X$ provient d'une $(\Bbb S^1)^m$-action ; c'est-\`a-dire, tenant compte de la condition (K0) si et seulement si l'action transversalement k\"ahl\'erienne sur $N$ est donn\'ee par une action de tores. Or d'apr\`es le th\'eor\`eme A de \cite{M-V}, c'est le cas si et seulement si la condition (K) est v\'erifi\'ee. Ceci montre le (ii).
$\square$
\enddemo

Notons pour finir le

\proclaim{Corollaire 20.5}
Sous les hypoth\`eses du th\'eor\`eme 20.2, si de plus $m$ vaut $1$, alors $X$ est sasakienne.
\endproclaim

\demo{Preuve}
En suivant de pr\`es la preuve du th\'eor\`eme 20.2, on voit que la forme transversalement k\"ahl\'erienne $\omega$ descend \`a $X$ ; et que sa projection est la dif\-f\'e\-ren\-ti\-el\-le d'une forme de contact.
$\square$
\enddemo

\remark{Remarque}
Contrairement au cas de la proposition 20.1, la vari\'et\'e LVM $N$ n'est pas en g\'en\'eral une vari\'et\'e associ\'ee au $\Bbb R^m$-c\^one de $X$, car le fibr\'e en tores $N\to X$ n'est pas un fibr\'e trivial.
\endremark

\remark{Remarque}
Ainsi que nous l'avons signal\' e plus haut, on peut facilement \'etablir une version orbifold du th\'eor\`eme 20.2, (i) en relaxant l'hypoth\`ese (ii) dans la d\'efinition de la condition (K0). Tous les r\'esultats de \cite{M-V} sont en effet \'enonc\'es pour le cas orbifold. Bien que plus technique, une telle version a l'avantage de correspondre \`a une condition (K0) plus facile \`a manier.
\endremark

\remark{Remarque}
Il n'est pas clair que le m\^eme type de r\'esultat soit vrai pour les vari\'et\'es de \cite{Bo}, qui g\'en\'eralisent les vari\'et\'es LVM. En effet, il n'y a pas d'\'equivalent connu \`a la forme $\omega$, car il n'y a pas forc\'ement de feuilletage {\it transversalement k\"ahl\'erien} sur un tel objet, cf. \cite{CF-Z}. 
\endremark

\vfill
\eject

\enddocument